\documentclass{amsart}
\usepackage[T1]{fontenc}
\usepackage[utf8]{inputenc}
\usepackage{amsmath}
\usepackage{amsfonts}
\usepackage{amssymb}
\usepackage{lmodern}
\usepackage{microtype}
\usepackage[hidelinks]{hyperref}
\hypersetup{pdftitle={Caffarelli--Kohn--Nirenberg inequalities in
rearrangement-invariant spaces},
  pdfauthor={Joaquim Martin}}
\numberwithin{equation}{section}

\newtheorem{theorem}{Theorem}[section]
\newtheorem{lemma}[theorem]{Lemma}
\newtheorem{proposition}[theorem]{Proposition}
\newtheorem{corollary}[theorem]{Corollary}

\theoremstyle{remark}
\newtheorem{remark}[theorem]{Remark}

\newcommand{\R}{\mathbb R}
\newcommand{\HH}{\mathcal H}
\newcommand{\supp}{\operatorname{supp}}
\newcommand{\Xbar}{\overline X}
\newcommand{\Ybar}{\overline Y}
\newcommand{\Zbar}{\overline Z}

\makeatletter \@namedef{subjclassname@2020}{\textup{2020}
Mathematics Subject Classification} \makeatother

\begin{document}

\title[CKN inequalities in r.i.\ spaces]
{Caffarelli--Kohn--Nirenberg inequalities in rearrangement-invariant
spaces}

\author{Joaquim Mart\'{\i}n}

\address{Department of Mathematics\\
Universitat Aut\`onoma de Barcelona\newline
Bellaterra (Barcelona), Spain}

\email{Joaquin.Martin@uab.cat}

\urladdr{https://orcid.org/0000-0002-7467-787X}

\author{Walter A.~Ortiz}

\address{Department of Mathematics\\
Universitat Aut\`onoma de Barcelona Bellaterra (Barcelona), Spain;
School of Engineering, Science and Technology, Universidad
Internacional de Valencia (VIU), Valencia, Spain; and Faculty of
Economic and Social Sciences, Universitat Internacional de Catalunya
(UIC Barcelona), Barcelona, Spain}

\email{walterandres.ortiz@professor.universidadviu.com}

\urladdr{https://orcid.org/0000-0002-8617-3919}



\begin{abstract}
We characterize the Caffarelli--Kohn--Nirenberg inequalities
\[
\big\||x|^{-b}f\big\|_Y \leq C\big\||x|^{-a}\nabla f\big\|_X, \qquad
f\in C_c^\infty(\mathbb R^n),\quad n\geq2,
\]
for arbitrary pairs of rearrangement-invariant Banach function
spaces $X,Y$ over Lebesgue measure, when $b\geq0$. The inequality
forces $a\leq b$ and, in that case, it holds if and only if an
explicit one-dimensional weighted Hardy operator, in which the two
weights remain distinct multipliers, is bounded between the
representation spaces of $X$ and $Y$. The results previously
available beyond the Lebesgue scale concern particular families of
spaces (Lorentz, weak Lebesgue, Orlicz).

When the gradient weight is singular ($a>0$), symmetrization with
respect to Lebesgue measure is not available, and the perturbative
reduction to an unweighted gradient requires a Hardy inequality in
$X$, which fails at the critical exponents. We overcome this
difficulty by combining relative rearrangements with an
integrating-factor argument, which yields a pointwise estimate
independent of the norms. Consequences include the characterization
of optimal rearrangement-invariant targets, critical embeddings with
optimal domain and range, Lorentz refinements, sharp weighted Hardy
inequalities, and a comparison between the optimal constants for
general and radial functions.
\end{abstract}

\maketitle

\section{Introduction}

In their seminal paper \cite{CKN}, Caffarelli, Kohn and Nirenberg
established a family of first-order interpolation inequalities with
power weights on $\R^n$, now known as the
Caffarelli--Kohn--Nirenberg (CKN) inequalities. The family contains
the Sobolev and Hardy inequalities as particular cases, and the
associated extremal problems, including the symmetry and symmetry
breaking of optimizers, have been studied extensively (see, for
example, \cite{CatrinaWang,FelliSchneider,DolbeaultEstebanLoss}). In
this paper we consider the Sobolev branch of the CKN inequalities,
that is, the two-norm inequalities
\begin{equation}\label{eq:intro-Lp}
 \big\||x|^{-b}f\big\|_{L^q(\R^n)}
 \leq C\big\||x|^{-a}\nabla f\big\|_{L^p(\R^n)},
 \qquad f\in C_c^\infty(\R^n),
\end{equation}
with $n\geq2$. Writing $\delta=b-a$, dilation invariance forces
$1/q=1/p-(1-\delta)/n$. The Sobolev inequality corresponds to
$a=b=0$, the Hardy inequality to $a=0$, $b=1$, and the
Hardy--Sobolev inequalities to $a=0$, $0<b<1$; the parameter
$\delta$ measures the balance between the two weights, from the
Sobolev regime $\delta=0$ to the Hardy regime $\delta=1$. In the
Lebesgue scale, the admissible parameters are completely described
in \cite{CKN}.

Analogues of these inequalities have since been obtained beyond the
Lebesgue scale, in particular families of spaces: Lorentz spaces
\cite{Yacoub}, Lorentz estimates derived from weighted Riesz
potentials \cite{NursultanovTikhonov,LiuYuZhou}, weak Lebesgue
spaces, including the interpolation form of the inequalities, by
sparse domination of fractional integrals \cite{Wang}, and modular
inequalities in Orlicz--Sobolev spaces \cite{BoceaMihailescu}; see
also \cite{FiorenzaEtAl} for Gagliardo--Nirenberg inequalities in
rearrangement-invariant spaces and \cite{DiFrattaFiorenza} for
radial weights. Each of these works is tied to a particular scale
and to tools specific to it (potential representations, sparse
domination, modular arguments); none of them provides a criterion
for general rearrangement-invariant spaces.

The purpose of this paper is to characterize the inequality
\begin{equation}\label{eq:intro}
 \big\||x|^{-b}f\big\|_Y
 \leq C\big\||x|^{-a}\nabla f\big\|_X,
 \qquad f\in C_c^\infty(\R^n),
\end{equation}
for an arbitrary pair of rearrangement-invariant Banach function
spaces\footnote{The basic notions and notation used in this
introduction and throughout the paper are recalled in
Section~\ref{sec:preliminaries}, or introduced where needed.}, that
is, to give a single criterion which is necessary and sufficient for
every such pair $X,Y$ over Lebesgue measure (r.i.\ spaces, for
short).

In (\ref{eq:intro}) the norms are invariant under rearrangements
with respect to Lebesgue measure, while the weights act as
multipliers. For an individual Lebesgue norm, multiplication by a
weight can be incorporated into the measure: for instance,
$\||x|^{-b}f\|_{L^q}$ is the norm of $f$ in $L^q(|x|^{-bq}dx)$. The
two sides of (\ref{eq:intro-Lp}) generally lead to different
measures, $|x|^{-bq}dx$ and $|x|^{-ap}dx$. Outside the Lebesgue
scale, even this identification of each multiplier norm with a norm
over a weighted measure is no longer available in general. Weighted
Sobolev inequalities have a long history, including two-weight
formulations \cite{Mazya}, inequalities in weighted Lebesgue spaces
\cite{CabreRosOton}, and Sobolev inequalities in r.i.\ spaces with
respect to a common weighted measure
\cite{MartinMilman,FeoMartinPosteraro,Drazny} and its references.

Inequality (\ref{eq:intro}) forces $a\leq b$
(Proposition~\ref{prop:a-le-b}). Our main result,
Theorem~\ref{thm:main}, shows that, for $b\geq0$ and $a\leq b$, it
holds if and only if the weighted Hardy operator
\[
\HH_{a,b}h(t) =t^{-b/n}\int_t^\infty h(s)s^{(a+1)/n-1}\,ds
\]
is bounded from $\Xbar$ to $\Ybar$, where $\Xbar,\Ybar$ are the
representation spaces of $X,Y$.

The two weights remain separate: the gradient weight determines the
kernel, while the function weight appears as an exterior factor.
Necessity follows by testing on radial functions. The main
difficulty is proving sufficiency for general functions,
particularly when $a>0$, where the standard weighted isoperimetric
approach does not apply directly.

A natural starting point for such a reduction is symmetrization. In
the unweighted setting, this method reduces Sobolev inequalities in
r.i.\ spaces to one-dimensional inequalities (see, for example,
\cite{MartinMilmanPustylnik,MartinMilman,ACPS,MihulaPRSE}). When two
different power weights are present, one may instead incorporate the
function weight into a reference measure and the gradient weight
into a corresponding weighted perimeter, and then use rearrangements
with respect to this measure together with weighted isoperimetric
inequalities (see \cite{AlvinoEtAl} and the references therein).
This approach, however, does not directly provide the reduction we
seek, since the norms in (\ref{eq:intro}) are taken with respect to
Lebesgue measure and the two weights must remain as separate
multipliers.

With Lebesgue rearrangements, the function weight causes no
difficulty when $b\geq0$: the Hardy--Littlewood inequality bounds
$\||x|^{-b}f\|_Y$ by the norm of $\omega_n^{b/n}t^{-b/n}f^*(t)$ in
$\Ybar$, with equality for radially non-increasing $f$. The gradient
weight is the obstacle. When $a\leq0$, an isoperimetric inequality
with the weight $|x|^{-a}$ on the perimeter and Lebesgue measure as
volume yields the reduction, even with the optimal radial constant
(see \cite{AlvinoEtAl} and Section~\ref{sec:preliminaries}). When
$a>0$, no such inequality exists: balls of fixed volume have
$|x|^{-a}$-weighted perimeter tending to zero as their centers move
to infinity. To our knowledge, a necessary-and-sufficient reduction
for arbitrary r.i.\ pairs, with both weights acting as multipliers
over Lebesgue measure, has not been available in this regime.

For $a>0$ there is, nevertheless, a simple perturbative route which
uses only rearrangements with respect to Lebesgue measure. Putting
$u=|x|^{-a}|f|$, one has $|\nabla u|\leq|x|^{-a}|\nabla
f|+a|x|^{-a-1}|f|$. Hence, if the weighted Hardy inequality
\begin{equation}\label{eq:intro-hardy}
\big\||x|^{-a-1}f\big\|_X\leq C\big\||x|^{-a}\nabla f\big\|_X
\end{equation}
is available in $X$, the problem reduces to the case of an
unweighted gradient (Proposition~\ref{prop:perturbative}). In the
Lebesgue scale this covers the whole subcritical range $a<n/p-1$.
However, (\ref{eq:intro-hardy}) is an additional hypothesis on $X$,
namely the boundedness of $\HH_{a,a+1}$ on $\Xbar$, and it is not
implied by the boundedness of $\HH_{a,b}:\Xbar\to\Ybar$. At the
critical exponent $a=n/p-1>0$, inequality (\ref{eq:intro}) holds
from $X=L^{p,1}$ into $Y=L^{q,\infty}$, whereas
(\ref{eq:intro-hardy}) fails for $X=L^{p,1}$
(Remark~\ref{rem:perturbative}). Thus the perturbative argument
recovers individual inequalities but cannot produce a
characterization, and a different method is needed precisely at such
endpoints.

To overcome this obstruction without any assumption on the spaces,
we combine a transfer of weights with
relative rearrangements, developed by Rakotoson and collaborators
(see \cite{RakotosonTemam,RakotosonPointwise,Rakotoson}). More
precisely, for $a>0$ we choose $a\leq c\leq\min\{b,a+1\}$ and
introduce $u=|x|^{-c}|f|$. Writing $\kappa=c-a$ and
$g=|x|^{-a}|\nabla f|$, a direct computation gives
\begin{equation}\label{eq:weight-transfer}
|x|^\kappa|\nabla u| \leq g+c|x|^{\kappa-1}u.
\end{equation}

The perimeter weight $|x|^\kappa$ is now admissible for isoperimetry
with Lebesgue volume, but the additional term $c|x|^{\kappa-1}u$
must be controlled while preserving the critical powers. Relative
rearrangements, denoted by $(\cdot)_{*u}$, provide this control: for
almost every $s$ such that $(u^*)'(s)<0$, we have
\[
(|x|^{\kappa-1}u)_{*u}(s) =u^*(s)\,(|x|^{\kappa-1})_{*u}(s)
\]
(Lemma~\ref{lem:active}). This factorization makes its contribution
linear in $u^*$, allowing the resulting differential inequality to
be integrated. A logarithmic estimate then controls the integrating
factor without changing the critical powers
(Lemma~\ref{lem:pointwise}).

The outcome is a pointwise bound for $u^*$ by a Hardy operator
applied to the relative rearrangement of $g$. Since relative
rearrangements are contractions in every r.i.\ space, this bound is
independent of the norms, which makes a characterization for
arbitrary pairs possible. In this sense, relative rearrangements
provide a way of symmetrizing two-weight inequalities without
changing the reference measure.

The method yields several further results.
Proposition~\ref{prop:sharp} compares the optimal constants for
general and radial functions: the former lies between the latter and
twice its value, and they coincide when $a\leq0$ or $\delta\geq1$.
The reduction also covers $\delta>1$, although admissible
pairs can exist only for $\delta<n+1$
(Remark~\ref{rem:effective-range}). For $X=L^p$ and $\delta>1$, the
target $L^q$ is never admissible, whereas $L^{q,p}$ is admissible
when $a<n/p-1$ and $q>1$. This boundedness is elementary (Hardy's
inequality followed by O'Neil's H\"older inequality); what the
reduction adds is that $L^{q,p}$ is the optimal
rearrangement-invariant target and that no such target exists when
$a\geq n/p-1$ (Proposition~\ref{cor:optimal-lebesgue} and
Remark~\ref{rem:delta-large}). More generally, Theorem~\ref{thm:optimal} characterizes
the existence of an r.i.\ target and identifies the optimal one,
while the theory of weighted Hardy operators on r.i.\ spaces
\cite{Mihula} transfers directly to CKN inequalities, including
those in Orlicz and Lorentz--Zygmund spaces.

Weighted Riesz potential estimates provide another route to CKN
inequalities in Lorentz and weak Lebesgue spaces, including some
weak-type endpoints \cite{NursultanovTikhonov,LiuYuZhou,Wang}. The present reduction also
covers the range $a<-n/p'$, beyond these potential estimates (see
Remark~\ref{rem:delta-large}), and treats $p=1$, with explicit
constants, and $q=\infty$. At the critical exponent $a=n/p-1$, it
identifies $L^{p,1}$ and $L^{q,\infty}$ as the optimal domain and
target among all r.i.\ spaces (Proposition~\ref{cor:critical});
there, no r.i.\ target exists for $X=L^p$ when $p>1$. In particular,
the weak-type Hardy inequality $\|f/|x|\|_{L^{n,\infty}}\leq
C\|\nabla f\|_{L^n}$ fails for $n\geq2$ (Remark~\ref{rem:critical}).

Throughout, we assume $b\geq0$. In the reduction argument, this
restriction enters through the Hardy--Littlewood step: for $b<0$ the
function weight is radially increasing, and $\||x|^{-b}f\|_Y$ is not
controlled by $f^*$ alone. In the Lebesgue scale, the classical
theorem \cite{CKN} still gives (\ref{eq:intro-Lp}) for $b<0$, $a\leq
b$ and $1\leq p\leq q<\infty$, provided $1/q=1/p-(1-\delta)/n$; the
remaining Hardy condition $1/q>b/n$ is then automatic
(Remark~\ref{rem:tests}). For $p=2$, however, symmetry breaking
occurs in a region where $b<0$, and the optimal constant for general
functions is strictly larger than the radial one
\cite{CatrinaWang,FelliSchneider,DolbeaultEstebanLoss}. Whether
Theorem~\ref{thm:main} extends to $b<0$ for arbitrary r.i.\ pairs
remains open.

Section~\ref{sec:preliminaries} collects the necessary facts on
r.i.\ spaces, relative rearrangements and weighted isoperimetry.
Section~\ref{sec:main} contains the pointwise estimate, the main
theorem, the perturbative argument and its limitations, and the
comparison of constants. Section~\ref{sec:optimal}
characterizes optimal targets, and Section~\ref{sec:examples} treats
Lebesgue and Lorentz spaces, critical embeddings and weighted Hardy
inequalities.

\section{Preliminaries}
\label{sec:preliminaries}

Throughout the paper, $A\preceq B$ means that $A\leq CB$ for a
constant $C>0$ which may depend on the fixed structural parameters
and spaces, but not on the functions under consideration. We write
$A\simeq B$ when both $A\preceq B$ and $B\preceq A$ hold.

\subsection{Rearrangements and rearrangement-invariant spaces}

Let $\Omega\subset\mathbb R^n$ be a measurable set, $n\geq2$,
equipped with Lebesgue measure, and put $0<m=|\Omega|\leq\infty$.
For a measurable function $f$ on $\Omega$, its distribution function
is defined by
\[
\mu_f(\lambda) =|\{x\in\Omega:|f(x)|>\lambda\}|, \qquad \lambda>0.
\]
Two measurable functions, possibly defined on different measure
spaces, are said to be equimeasurable if their distribution
functions coincide.

The decreasing rearrangement and the maximal rearrangement of $f$
are defined, respectively, by
\[
f^*(t)=\inf\{\lambda\geq0:\mu_f(\lambda)\leq t\}, \qquad
f^{**}(t)=\frac1t\int_0^t f^*(s)\,ds, \qquad 0<t<m.
\]
When $m=\infty$, the interval $(0,m)$ is understood as $(0,\infty)$.

A Banach function space $X=X(\Omega)$ is understood to be equipped
with a Banach function norm in the sense of \cite[Chapter~1]{BS}. In
particular, its norm satisfies the Fatou property, is finite on
characteristic functions of sets of finite measure, and controls
integration over such sets. We say that $X$ is
rearrangement-invariant (r.i.) if, whenever $f\in X$ and $g$ is
measurable on $\Omega$, then
\begin{equation*}
f^*=g^* \quad\Longrightarrow\quad g\in X,\qquad \|g\|_X=\|f\|_X.
\end{equation*}

Classical examples of r.i.\ spaces include Lebesgue, Lorentz,
Lorentz--Zygmund and Orlicz spaces. Every r.i.\ space $X=X(\Omega)$
admits a representation space $\Xbar$ on $(0,m)$, where
$m=|\Omega|$, such that
\begin{equation}\label{eq:representation}
\|f\|_X=\|f^*\|_{\Xbar}.
\end{equation}
Conversely, every r.i.\ Banach function space $\Xbar$ on $(0,m)$
determines an r.i.\ Banach function space $X$ on $\Omega$ through
(\ref{eq:representation}); see \cite[Chapter~2, Theorem~4.10 and
subsequent remarks]{BS}.

For a Banach function space $X$ on $\Omega$, its associate space
(K\"othe dual) $X'$ is defined by
\begin{equation*}
\|g\|_{X'} =\sup_{\|f\|_X\leq1} \int_\Omega |f(x)g(x)|\,dx.
\end{equation*}
If $X$ is rearrangement-invariant, the Hardy--Littlewood inequality
and the rearrangement principle for associate norms give
\[
\begin{aligned}
\|g\|_{X'} &=\sup_{\|f\|_X\leq1}
  \int_\Omega |f(x)g(x)|\,dx\\
&=\sup_{\|\varphi\|_{\Xbar}\leq1}
  \int_0^m |\varphi(t)|g^*(t)\,dt\\
&=\|g^*\|_{(\Xbar)'}.
\end{aligned}
\]
Thus $X'$ is also rearrangement-invariant, and its representation
space satisfies
\begin{equation*}
\overline{X'}=(\Xbar)'
\end{equation*}
isometrically. Moreover, the Lorentz--Luxemburg theorem gives
$X''=X$ isometrically, without requiring absolute continuity of the
norm; see \cite[Chapters~1 and 2]{BS}.

For nonnegative measurable functions $h,k$ on $\Omega$, we write
$h\prec\!\prec k$ if
\begin{equation*}
\int_0^t h^*(s)\,ds \leq\int_0^t k^*(s)\,ds, \qquad 0<t<m.
\end{equation*}
The same notation is used for measurable functions on $(0,m)$. By
the Hardy--Littlewood--P\'olya principle (see \cite[Chapter~2,
Theorem~4.6]{BS}),
\begin{equation}\label{eq:HLP}
h\prec\!\prec k \quad\Longrightarrow\quad \|h\|_E\leq\|k\|_E
\end{equation}
for every r.i.\ Banach function space $E$ on the underlying measure
space, provided $h$ and $k$ are defined on that same space. In view of
(\ref{eq:representation}), the same conclusion holds, in the form
$\|h\|_X\leq\|k\|_{\Xbar}$, when $h$ is defined on $\Omega$ and $k$
on $(0,m)$; we shall use (\ref{eq:HLP}) in this sense without
further comment.

When $X=X(\mathbb R^n)$ and $\Omega\subset\mathbb R^n$ has finite
measure, we denote by $X(\Omega)$ the space induced by $X$, equipped
with the norm
\[
\|h\|_{X(\Omega)} =\|\widetilde h\|_{X(\mathbb R^n)},
\]
where $\widetilde h$ is the extension of $h$ by zero outside
$\Omega$.

We recall the Lorentz spaces used below (see \cite[Chapter~4,
Section~4]{BS}).  For $1<p<\infty$ and $1\leq r<\infty$, the space
$L^{p,r}(\Omega)$ consists of all measurable functions $f$ such that
\[
\|f\|_{L^{p,r}}^* = \left( \int_0^m
\bigl[t^{1/p}f^*(t)\bigr]^r\frac{dt}{t} \right)^{1/r} <\infty.
\]
For $r=\infty$, we set $\|f\|_{L^{p,\infty}}^*
=\sup_{0<t<m}t^{1/p}f^*(t)$. When $1<p<\infty$, replacing
$f^*$ by $f^{**}$ gives equivalent norms. We denote these
$f^{**}$-based norms by $\|f\|_{L^{p,r}}$ for $1<r\leq\infty$;
for $r=1$ we use the $f^*$-based normalization in
(\ref{eq:Lp1}). In particular, for $q>1$,
\begin{equation}
\label{eq:weak-norm} \|f\|_{L^{q,\infty}}^* \leq
\|f\|_{L^{q,\infty}} \leq q'\|f\|_{L^{q,\infty}}^*, \qquad
\|f\|_{L^{q,\infty}} =\sup_{0<t<m}t^{1/q}f^{**}(t).
\end{equation}

For the critical results below, we use the normalization
\begin{equation}
\label{eq:Lp1} \|f\|_{L^{p,1}} =\int_0^m f^*(t)t^{1/p-1}\,dt, \qquad
1\leq p<\infty.
\end{equation}
In particular, $L^{1,1}=L^1$. The equivalence with the expressions
involving $f^{**}$ stated above is used only for $p>1$.

For further background on rearrangements and r.i.\ spaces, see
\cite{BS,KPS}.

\subsection{Relative rearrangements}

We recall the definition and properties of relative rearrangements
needed below; see \cite[Chapter~2]{Rakotoson} and
\cite{RakotosonPointwise,RakotosonSimon}. Their connection with
coarea and pointwise Sobolev estimates is developed in
\cite{RakotosonTemam}.

Let $\Omega\subset\mathbb R^n$ be a bounded open set, let $u\in
L^1(\Omega)$ be nonnegative, and put $m=|\Omega|$. For $h\in
L^1(\Omega)$, $h\geq0$, and $0<s<m$, define
\begin{equation}
\label{eq:Rh} R_h(s) =\int_{\{u>u^*(s)\}}h(x)\,dx
+\int_0^{s-|\{u>u^*(s)\}|} \bigl(h|_{\{u=u^*(s)\}}\bigr)^*(t)\,dt.
\end{equation}
The rearrangement in the second integral is taken with respect to
Lebesgue measure restricted to the level set $\{u=u^*(s)\}$.

The function $R_h$ is absolutely continuous on $[0,m]$, and its
almost everywhere derivative
\[
h_{*u}(s)=R_h'(s)
\]
is called the relative rearrangement of $h$ with respect to $u$. In
particular,
\begin{equation*}
 R_h(s)=\int_0^s h_{*u}(t)\,dt, \qquad 0\leq
s\leq m.
\end{equation*}
 The representation (\ref{eq:Rh}) and the absolute
continuity of $R_h$ follow from the directional derivative formula
for rearrangements; see \cite[Theorems~2.1.2 and 2.1.3 and
Definition~2.1.3]{Rakotoson}.

The second term in (\ref{eq:Rh}) accounts for level sets of positive
measure. In particular, if $|\{u=u^*(s)\}|=0$, this term vanishes,
and hence
\begin{equation}
\label{eq:Rh-active} R_h(s)=\int_{\{u>u^*(s)\}}h(x)\,dx.
\end{equation}

The relative rearrangement is order-preserving: $h\leq k$ implies
$h_{*u}\leq k_{*u}$ almost everywhere on $(0,m)$. It also satisfies
the submajorization property
\begin{equation}
\label{eq:rel-submaj} \int_0^t (h_{*u})^*(s)\,ds \leq \int_0^t
h^*(s)\,ds, \qquad 0<t<m.
\end{equation}
Both properties follow from \cite[Proposition~2.2.1(b) and
Corollary~2.2.1(1)]{Rakotoson}. Consequently, extending $h_{*u}$ by
zero outside $(0,m)$, the Hardy--Littlewood--P\'olya principle
(\ref{eq:HLP}) yields
\begin{equation}
\label{eq:rel-contract} \|h_{*u}\|_{\overline X} \leq
\|h\|_{X(\Omega)}
\end{equation}
for every r.i.\ Banach function space $X=X(\mathbb R^n)$.

Although the relative rearrangement is generally nonlinear, the
following lemma will show that useful additivity and product
identities can be recovered almost everywhere on the set of active
levels, where $(u^*)'<0$.

\begin{lemma}\label{lem:active}
Let $\Omega\subset\mathbb R^n$ be a bounded open set and let
$u\geq0$ belong to $W^{1,1}_0(\Omega)$. Put $m=|\Omega|$ and
\[
J=\{s\in(0,m):(u^*)'(s)\ \text{exists and}\ (u^*)'(s)<0\}.
\]
Then, for every $s\in J$,
\[
|\{u=u^*(s)\}|=0, \qquad |\{u>u^*(s)\}|=s.
\]

Moreover, if $h,k\in L^1(\Omega)$ are nonnegative and $uh\in
L^1(\Omega)$, then, for almost every $s\in J$,
\begin{equation}\label{eq:rel-identities}
(h+k)_{*u}(s)=h_{*u}(s)+k_{*u}(s), \qquad
(uh)_{*u}(s)=u^*(s)h_{*u}(s).
\end{equation}
\end{lemma}

\begin{proof}

Since $u\in W^{1,1}_0(\Omega)$, its decreasing rearrangement $u^*$
is locally absolutely continuous on $(0,m)$; see, for example,
\cite[Corollary~3.3.1]{Rakotoson} and \cite{Talenti,BrothersZiemer}.
As $u^*$ is non-increasing, $(u^*)'\leq0$ almost everywhere, with
$(u^*)'=0$ almost everywhere on $(0,m)\setminus J$. It follows that
\[
u^*(s)-u^*(t) =-\int_{(s,t)\cap J}(u^*)'(\tau)\,d\tau, \qquad
0<s<t<m.
\]
Thus all the variation of $u^*$ is accounted for by $J$.

We first identify the levels corresponding to points of $J$. By the
definition of the decreasing rearrangement,
\begin{equation}\label{eq:level-bounds}
|\{u>u^*(s)\}| \leq s\leq |\{u\geq u^*(s)\}|, \qquad 0<s<m.
\end{equation}

Fix $s\in J$ and put $\lambda=u^*(s)$. Suppose that
$|\{u=\lambda\}|>0$. Then
\begin{equation*}
I_\lambda= \bigl(|\{u>\lambda\}|,|\{u\geq\lambda\}|\bigr)
\end{equation*}
is a nonempty interval on which $u^*$ is constant and equal to
$\lambda$. By (\ref{eq:level-bounds}), $s\in\overline{I_\lambda}$,
and since $u^*$ is differentiable at $s$, its constancy on
$I_\lambda$ forces $(u^*)'(s)=0$, contradicting $s\in J$. Therefore,
for every $s\in J$,
\[
|\{u=u^*(s)\}|=0, \qquad |\{u>u^*(s)\}|=s.
\]
In particular, for every nonnegative $\ell\in L^1(\Omega)$, the
relative primitive $R_\ell(s)$ has the representation
(\ref{eq:Rh-active}) at every $s\in J$.

Moreover, every $s\in J$ is a right accumulation point of $J$.
Indeed, since $(u^*)'(s)<0$, there exists $\varepsilon_0>0$,
depending on $s$, such that $s+\varepsilon_0<m$ and
$u^*(s+\varepsilon)<u^*(s)$ whenever $0<\varepsilon<\varepsilon_0$.
By the absolute continuity of $u^*$ on $[s,s+\varepsilon]$, we
obtain
\[
0<u^*(s)-u^*(s+\varepsilon) =-\int_{(s,s+\varepsilon)\cap
J}(u^*)'(\tau)\,d\tau.
\]
Thus $|J\cap(s,s+\varepsilon)|>0$ for every
$0<\varepsilon<\varepsilon_0$. In particular, there exists a
sequence $s_j\in J$ such that $s_j>s$ and $s_j\to s$.

We now establish the first identity in (\ref{eq:rel-identities}). By
(\ref{eq:Rh-active}) and the linearity of the integral,
\[
R_{h+k}(s)=R_h(s)+R_k(s), \qquad s\in J,
\]
as all three integrals are taken over the same set $\{u>u^*(s)\}$.

Set
\[
F=R_{h+k}-R_h-R_k.
\]
The integrability of $h$ and $k$ ensures that the relative
primitives, and hence $F$, are absolutely continuous. Moreover,
$F=0$ on $J$.

Let $s\in J$ be a point where $F$ is differentiable. The
accumulation property established above provides a sequence $s_j\in
J$ with $s_j>s$ and $s_j\to s$. As $F(s_j)=F(s)=0$,
\[
F'(s)=\lim_{j\to\infty} \frac{F(s_j)-F(s)}{s_j-s}=0.
\]
The absolute continuity of $F$ guarantees differentiability almost
everywhere. It follows that
\[
R_{h+k}'(s)=R_h'(s)+R_k'(s) \qquad\text{for almost every }s\in J.
\]
Recalling that $R_\ell'=\ell_{*u}$ almost everywhere, we obtain
\[
(h+k)_{*u}(s)=h_{*u}(s)+k_{*u}(s) \qquad\text{for almost every }s\in
J.
\]
This proves the first identity in (\ref{eq:rel-identities}).

For the second identity, let $s\in J$ be a point where $R_h$ and
$R_{uh}$ are differentiable. Their absolute continuity ensures that
this holds for almost every $s\in J$. By the accumulation property
established above, we may approach $s$ from the right through points
$t\in J$. For such $t$, the representation (\ref{eq:Rh-active}) and
the monotonicity of $u^*$ give
\[
\begin{aligned}
R_h(t)-R_h(s)
&=\int_{\{u^*(t)<u\leq u^*(s)\}}h(x)\,dx,\\
R_{uh}(t)-R_{uh}(s) &=\int_{\{u^*(t)<u\leq u^*(s)\}}u(x)h(x)\,dx.
\end{aligned}
\]
On the common domain of integration we have $u^*(t)\leq u(x)\leq
u^*(s)$, and therefore, as $h\geq0$,
\[
\begin{aligned}
u^*(t)\bigl(R_h(t)-R_h(s)\bigr)
&\leq R_{uh}(t)-R_{uh}(s)\\
&\leq u^*(s)\bigl(R_h(t)-R_h(s)\bigr).
\end{aligned}
\]
Dividing by $t-s$ and letting $t\downarrow s$ through points of $J$,
the continuity of $u^*$ and the differentiability of $R_h$ and
$R_{uh}$ yield
\[
R_{uh}'(s)=u^*(s)R_h'(s).
\]
Recalling that $R_\ell'=\ell_{*u}$ almost everywhere, we obtain
\[
(uh)_{*u}(s)=u^*(s)h_{*u}(s) \qquad\text{for almost every }s\in J,
\]
which proves the second identity in (\ref{eq:rel-identities}).
\end{proof}

\begin{remark}\label{rem:product}
The product identity in (\ref{eq:rel-identities}) is not restricted
to active levels. For nonnegative $u,h\in L^1(\Omega)$ with $uh\in
L^1(\Omega)$ and $0<s<t<m$, the representation (\ref{eq:Rh}),
together with $\bigl((uh)|_{\{u=\lambda\}}\bigr)^*
=\lambda\bigl(h|_{\{u=\lambda\}}\bigr)^*$, gives
\[
u^*(t)\bigl(R_h(t)-R_h(s)\bigr)\leq R_{uh}(t)-R_{uh}(s)\leq
u^*(s)\bigl(R_h(t)-R_h(s)\bigr),
\]
and the right-continuity of $u^*$ yields
$(uh)_{*u}=u^*h_{*u}$ almost everywhere on $(0,m)$. Additivity, in
contrast, may fail on level sets of positive measure, and this is
the reason for restricting Lemma~\ref{lem:active} to $J$.
\end{remark}

\subsection{Weighted perimeter and isoperimetry}

We recall the notion of weighted perimeter needed for the relative
rearrangement--coarea identity and the corresponding weighted
isoperimetric inequality.

Let $\kappa\geq0$. For a bounded measurable set
$E\subset\mathbb R^n$, we define its weighted perimeter by
\begin{equation}\label{eq:perimeter}
P_\kappa(E) = \sup\left\{
\int_E\operatorname{div}(|x|^\kappa\Phi(x))\,dx: \Phi\in
C_c^1(\mathbb R^n;\mathbb R^n),\ |\Phi|\leq1 \right\}.
\end{equation}

For a bounded open set $\Omega\subset\mathbb R^n$, the relative
weighted perimeter $P_{\Omega,\kappa}(E)$, $E\subset\Omega$, is
defined by restricting the test vector fields in
(\ref{eq:perimeter}) to $C_c^1(\Omega;\mathbb R^n)$ (see
\cite[Definition~7.1.1]{Rakotoson}). Since these fields vanish near
$\partial\Omega$, the relative perimeter measures only the weighted
boundary of $E$ inside $\Omega$, without counting any contribution
on $\partial\Omega$. In particular, if $E$ is contained, up to a
null set, in a compact subset of $\Omega$, then
\begin{equation}\label{eq:perimeter-local}
P_{\Omega,\kappa}(E)=P_\kappa(E).
\end{equation}

The weighted coarea formula (see \cite[Lemma~7.1.1]{Rakotoson})
leads to the following relative rearrangement identity (see
\cite[Theorem~7.1.1]{Rakotoson}): for every nonnegative $u\in
W^{1,1}(\Omega)$ with $u^*$ locally absolutely
continuous\footnote{This condition is automatic for $u\in
W^{1,1}_0(\Omega)$; see \cite[Corollary~3.3.1]{Rakotoson}. In
particular, it holds in all our applications.},
\begin{equation}\label{eq:relative-coarea}
\bigl(|x|^\kappa|\nabla u|\bigr)_{*u}(s) =-(u^*)'(s)\,
P_{\Omega,\kappa}(\{u>u^*(s)\})
\end{equation}
for almost every $s$ such that $(u^*)'(s)<0$.

In our applications, $\Omega$ is a ball containing the compact
support of $u$ in its interior. The weight $|x|^\kappa$ belongs to
$W^{1,1}(\Omega)\cap L^\infty(\Omega)$, and the relative weighted
perimeter of every positive superlevel set coincides with its global
weighted perimeter.

The weighted isoperimetric inequality (see
\cite[Theorem~1.1]{AlvinoEtAl}), with $k=\kappa$ and $l=0$, gives
\begin{equation}\label{eq:weighted-iso}
P_\kappa(E) \geq n\omega_n^{(1-\kappa)/n}|E|^{(n-1+\kappa)/n},
\qquad \kappa\geq0,
\end{equation}
for every bounded measurable set $E\subset\mathbb R^n$, where
$\omega_n$ denotes the volume of the unit ball. Equality holds when
$E$ is a ball centered at the origin.

For $\kappa=0$, (\ref{eq:weighted-iso}) is the classical
isoperimetric inequality. For $\kappa\geq1$, it follows from the
elementary argument given below (for $\kappa=1$, see also
\cite[Theorem~5.1]{AlvinoEtAl}). For $0<\kappa<1$, the inequality is
case~(iii) of \cite[Theorem~1.1]{AlvinoEtAl} when $n\geq3$ (see
\cite[Remark~5.2]{AlvinoEtAl}), and case~(iv) when $n=2$, which
relies on a result of Csat\'o (see \cite{Csato} and \cite[Lemma~5.4
and Remark~5.3]{AlvinoEtAl}).

Although \cite[Theorem~1.1]{AlvinoEtAl} is stated for smooth sets,
the inequality extends to all measurable sets of finite measure,
with the distributional weighted perimeter; see \cite[Section~3,
particularly (3.7)]{AlvinoEtAl}. The approximation by smooth sets
in that reference gives convergence of both Lebesgue measure and
weighted perimeter.

For $\kappa\geq1$, let $E$ be a bounded measurable set, let
$\psi\in C_c^\infty(\mathbb R^n)$ satisfy $0\leq\psi\leq1$ and
$\psi=1$ on a neighbourhood of $\overline E$, and put
$\Phi_\varepsilon(x)=\psi(x)\,x\,(|x|^2+\varepsilon^2)^{-1/2}$ for
$\varepsilon>0$. Then $|\Phi_\varepsilon|\leq1$ and, on $E$,
\[
\operatorname{div}\bigl(|x|^\kappa\Phi_\varepsilon\bigr)
=\frac{(n+\kappa)|x|^\kappa}{(|x|^2+\varepsilon^2)^{1/2}}
-\frac{|x|^{\kappa+2}}{(|x|^2+\varepsilon^2)^{3/2}}
\geq\frac{(n+\kappa-1)|x|^\kappa}{(|x|^2+\varepsilon^2)^{1/2}}.
\]
Letting $\varepsilon\downarrow0$ and using monotone convergence, we
obtain
\begin{align*}
P_\kappa(E)
&\geq(n+\kappa-1)\int_E|x|^{\kappa-1}\,dx\\
&\geq(n+\kappa-1)\int_{E^\star}|x|^{\kappa-1}\,dx\\
&=n\omega_n^{(1-\kappa)/n}|E|^{(n-1+\kappa)/n},
\end{align*}
where $E^\star$ is the centered ball with $|E^\star|=|E|$; the second
inequality holds because $|x|^{\kappa-1}$ is radially
non-decreasing.

\section{The main reduction theorem}
\label{sec:main}

In this section we prove our main reduction theorem, which
characterizes the Caffarelli--Kohn--Nirenberg inequality in terms of
the boundedness of a Hardy operator. We begin with a pointwise
rearrangement estimate for inequalities of the form
(\ref{eq:weight-transfer}). The key step is to control the
zeroth-order term by combining the relative rearrangement identities
of Lemma~\ref{lem:active} with an integrating-factor argument.

For real numbers $a$ and $b$, we write
\begin{equation}\label{eq:Hardy-operator}
\HH_{a,b}h(s) =s^{-b/n}\int_s^\infty h(t)t^{(a+1)/n-1}\,dt, \qquad
s>0,
\end{equation}
for nonnegative measurable $h$ on $(0,\infty)$.

For $0<\gamma\leq1$, we also set
\begin{equation}
\label{eq:Theta} \Theta(\gamma)
=\int_0^1\frac{x^{\gamma-1}-1}{1-x}\,dx =\sum_{k=0}^\infty
\left(\frac1{k+\gamma}-\frac1{k+1}\right),
\end{equation}
where the series follows by expanding $(1-x)^{-1}$ and monotone
convergence. In terms of the digamma function $\psi$,
$\Theta(\gamma)=\psi(1)-\psi(\gamma)$.

The integrand is nonnegative and non-increasing in $\gamma$, so
$\Theta$ is non-increasing and $\Theta(1)=0$. Moreover, since
$x^\gamma<1$ for $0<x<1$,
\[
\frac{x^{\gamma-1}-1}{1-x}<x^{\gamma-1},
\]
and hence
\begin{equation}
\label{eq:Theta-bound} 0\leq\Theta(\gamma) <\int_0^1x^{\gamma-1}\,dx
=\frac1\gamma, \qquad 0<\gamma\leq1.
\end{equation}

\begin{lemma}
\label{lem:pointwise} Let $n\geq2$, $c\geq0$ and $\kappa\geq0$,
with $\kappa\leq1$ if $c>0$, and put $a=c-\kappa$ and
$\gamma=(n-1+\kappa)/n$. Let $u$ be a
nonnegative compactly supported function in $W^{1,1}(\R^n)$. Suppose
that $g\geq0$ is locally integrable and
\begin{equation}
\label{eq:hyp} |x|^\kappa|\nabla u| \leq g+c|x|^{\kappa-1}u
\quad\text{a.e. in }\R^n.
\end{equation}
Let $B$ be a ball centered at the origin containing the support of
$u$ in its interior, and let $g_{*u}$ denote the relative
rearrangement of $g$ computed in $B$, extended by zero beyond $|B|$.
Then
\begin{equation}
\label{eq:pointwise} u^*(s)\leq C_{n,\kappa,c}\,
\HH_{a,c}(g_{*u})(s) \quad\text{for every }s>0,
\end{equation}
where
\begin{equation}
\label{eq:constant} C_{n,\kappa,c}=
\frac{1}{n\omega_n^{(1-\kappa)/n}}
\begin{cases}
1, & c=0,\\
\exp(c\,\Theta(\gamma)/n), & c>0.
\end{cases}
\end{equation}
In particular, $C_{n,1,c}=1/n$ and $C_{n,\kappa,0}=1/\bigl(n
\omega_n^{(1-\kappa)/n}\bigr)$.
\end{lemma}

\begin{proof}
Put $m=|B|$ and $S=|\{u>0\}|$. The case $S=0$ is trivial, so we
assume $S>0$.

Since $u\in W^{1,1}_0(B)$, \cite[Corollary~3.3.1]{Rakotoson} gives
$u^*\in W^{1,1}(\varepsilon,m)$ for every $0<\varepsilon<m$.
Moreover, $S<m$ and $u^*=0$ on $(S,m)$. Hence $u^*$ is continuous at
$S$, and therefore
\begin{equation*}\label{eq:terminal}
 \lim_{t\uparrow S}u^*(t)=0.
\end{equation*}

We first check that the terms in (\ref{eq:hyp}) are integrable on
$B$. Since $B$ is bounded, $|x|^\kappa|\nabla u|$ and $g$ belong to
$L^1(B)$. If $\kappa>1$, then $|x|^{\kappa-1}$ is bounded on $B$,
and hence $|x|^{\kappa-1}u\in L^1(B)$. If $0<\kappa\leq1$, then
$|x|^{\kappa-1}\in L^n(B)$, and
the Sobolev embedding
\[
W^{1,1}(B)\hookrightarrow L^{n/(n-1)}(B)
\]
gives $|x|^{\kappa-1}u\in L^1(B)$.

For $\kappa=0$, the required integrability follows from the
classical Hardy inequality
\begin{equation}
\label{eq:hardy1} (n-1)\int_{\R^n}\frac{|\varphi(x)|}{|x|}\,dx \leq
\int_{\R^n}|\nabla\varphi(x)|\,dx, \qquad \varphi\in
C_c^\infty(\R^n),
\end{equation}
see \cite[Section~1.3.1, formula~(1.3.3)]{Mazya}. The extension of
(\ref{eq:hardy1}) to compactly supported $W^{1,1}(\R^n)$ functions
follows by standard approximation and Fatou's lemma. Hence
$|x|^{-1}u\in L^1(B)$.

Set
\[
J=\{s\in(0,S):(u^*)'(s)\ \text{exists and}\ (u^*)'(s)<0\}.
\]

Let $P_\kappa$ be the weighted perimeter defined in
(\ref{eq:perimeter}), and put
\[
A(s)=P_\kappa(\{u>u^*(s)\}).
\]

For $0<s<S$, the set $\{u>u^*(s)\}$ is contained, up to a null set,
in the compact support of $u$. Hence, by (\ref{eq:perimeter-local}),
\[
P_{B,\kappa}(\{u>u^*(s)\}) =P_\kappa(\{u>u^*(s)\})=A(s).
\]
Applying (\ref{eq:relative-coarea}) and (\ref{eq:perimeter-local}),
we obtain
\[
-A(s)(u^*)'(s) =\bigl(|x|^\kappa|\nabla u|\bigr)_{*u}(s)
\]
for almost every $s\in J$.

Set
\[
z=(|x|^{\kappa-1})_{*u}.
\]
It follows from the monotonicity of relative rearrangements and
Lemma~\ref{lem:active} that, for almost every $s\in J$,
\begin{equation}
\label{eq:diffA}
\begin{aligned}
-A(s)(u^*)'(s)
&=\bigl(|x|^\kappa|\nabla u|\bigr)_{*u}(s)\\
&\leq
\bigl(g+c|x|^{\kappa-1}u\bigr)_{*u}(s)\\
&=g_{*u}(s)+c\,u^*(s)z(s).
\end{aligned}
\end{equation}

Put
\[
A_0=n\omega_n^{(1-\kappa)/n}.
\]
By Lemma~\ref{lem:active}, we have $|\{u>u^*(s)\}|=s$ for every
$s\in J$. The weighted isoperimetric inequality
(\ref{eq:weighted-iso}) therefore gives
\begin{equation}\label{eq:A-lower}
A(s)\geq A_0s^\gamma, \qquad s\in J.
\end{equation}

Define
\[
\beta(s)=
\begin{cases}
\dfrac{z(s)}{A_0s^\gamma}, & s\in J,\\[6pt]
0, & s\in(0,S)\setminus J.
\end{cases}
\]
As $u^*$ is non-increasing, $(u^*)'=0$ almost everywhere on
$(0,S)\setminus J$. Combining (\ref{eq:diffA}) with
(\ref{eq:A-lower}), we obtain
\begin{equation}\label{eq:diff}
-(u^*)'(s) \leq A_0^{-1}s^{-\gamma}g_{*u}(s) +c\,\beta(s)u^*(s),
\qquad\text{for a.e. }s\in(0,S).
\end{equation}
Moreover, $\beta$ is measurable and, since $z\in L^1(0,m)$,
$\beta\in L^1_{\mathrm{loc}}(0,S)$.

Fix $0<s<t<S$ and define the integrating factor
\[
\mu(r)= \exp\left(c\int_s^r\beta(\tau)\,d\tau\right).
\]
Since $\mu'(r)=c\beta(r)\mu(r)$, (\ref{eq:diff}) gives
\[
\frac{d}{dr}\bigl(\mu(r)u^*(r)\bigr) \geq
-\frac1{A_0}\mu(r)r^{-\gamma}g_{*u}(r) \quad\text{for a.e.
}r\in(s,t).
\]

Since $\mu$ and $u^*$ are absolutely continuous on $[s,t]$,
integrating from $s$ to $t$ and using $\mu(s)=1$, we obtain
\begin{equation}
\label{eq:integrating-factor}
\begin{aligned}
u^*(s)\leq{}&
u^*(t)\exp\left(c\int_s^t\beta(\tau)\,d\tau\right)\\
&+\frac1{A_0}\int_s^t r^{-\gamma}g_{*u}(r)
\exp\left(c\int_s^r\beta(\tau)\,d\tau\right)\,dr.
\end{aligned}
\end{equation}

If $c=0$, the exponential factors in (\ref{eq:integrating-factor})
are identically one. If $c>0$, so that $0\leq\kappa\leq1$, we claim
that, for $0<s<r<S$,
\begin{equation}
\label{eq:beta-claim} \int_s^r\beta(\tau)\,d\tau \leq
\frac1n\log\frac rs+\frac{\Theta(\gamma)}n.
\end{equation}

Letting $t\uparrow S$ in (\ref{eq:integrating-factor}), we obtain
\begin{align*}
u^*(s) &\leq C_{n,\kappa,c} \int_s^S
r^{-\gamma}g_{*u}(r)
\left(\frac rs\right)^{c/n}\,dr\\
&=C_{n,\kappa,c}\,s^{-c/n} \int_s^S g_{*u}(r)r^{c/n-\gamma}\,dr.
\end{align*}
Since
\[
\frac cn-\gamma=\frac{a+1}{n}-1,
\]
the nonnegativity of $g_{*u}$ yields
\[
u^*(s)\leq C_{n,\kappa,c}\,\HH_{a,c}(g_{*u})(s), \qquad 0<s<S.
\]

For $s\geq S$, the same inequality is immediate because $u^*(s)=0$.
Thus (\ref{eq:pointwise}) follows once (\ref{eq:beta-claim}) is
established.

To prove (\ref{eq:beta-claim}) assume $c>0$, so that
$0\leq\kappa\leq1$.
 Put
$z_0=\bigl(|x|^{\kappa-1}\chi_B\bigr)^*$. If $\kappa<1$, then
$|\{x\in B:|x|^{\kappa-1}>\lambda\}|
\leq\omega_n\lambda^{-n/(1-\kappa)}$ for $\lambda>0$, while
$z_0=\chi_{(0,m)}$ if $\kappa=1$. In both cases,
\begin{equation}\label{eq:z0}
z_0(\sigma)\leq \omega_n^{(1-\kappa)/n}\sigma^{\gamma-1}, \qquad
0<\sigma<m.
\end{equation}
By definition, $0\leq\beta(\tau)\leq
A_0^{-1}\tau^{-\gamma}z(\tau)$ for almost every $\tau\in(0,S)$.

Fix $0<s<r<S$ and put $L=r-s$. The decreasing rearrangement of
$\tau\mapsto\tau^{-\gamma}\chi_{(s,r)}(\tau)$ is
$W(\sigma)=(s+\sigma)^{-\gamma}\chi_{(0,L)}(\sigma)$. Applying the
Hardy--Littlewood inequality (see, for example,
\cite[Chapter~2]{BS}), we obtain
\[
\int_s^r z(\tau)\tau^{-\gamma}\,d\tau \leq\int_0^\infty
F(\sigma)W(\sigma)\,d\sigma, \qquad F=(z\chi_{(s,r)})^*.
\]
Since $z\chi_{(s,r)}\leq z$, the submajorization
(\ref{eq:rel-submaj}) gives
\[
\int_0^t F(\sigma)\,d\sigma \leq\int_0^t z^*(\sigma)\,d\sigma
\leq\int_0^t z_0(\sigma)\,d\sigma, \qquad t>0.
\]
As $W$ is nonnegative and non-increasing, Hardy's lemma then yields
$\int_0^\infty FW\leq\int_0^\infty z_0W$. Combining these estimates
with (\ref{eq:z0}) and making the substitution $\sigma=sv$, we
obtain
\begin{align*}
\int_s^r\beta(\tau)\,d\tau &\leq\frac1{A_0}\int_0^L
z_0(\sigma)(s+\sigma)^{-\gamma}\,d\sigma\\
&\leq\frac1n\int_0^L
\sigma^{\gamma-1}(s+\sigma)^{-\gamma}\,d\sigma\\
&=\frac1n\int_0^{r/s-1} v^{\gamma-1}(1+v)^{-\gamma}\,dv.
\end{align*}

To estimate the last integral, observe that
\[
v^{\gamma-1}(1+v)^{-\gamma} =\frac1{1+v}
\left(\frac{v}{1+v}\right)^{\gamma-1} \geq\frac1{1+v}.
\]
Moreover, the substitution $x=v/(1+v)$ gives
\[
\begin{aligned}
&\int_0^\infty \left[ v^{\gamma-1}(1+v)^{-\gamma} -\frac1{1+v}
\right]\,dv\\
&\qquad= \int_0^1\frac{x^{\gamma-1}-1}{1-x}\,dx =\Theta(\gamma).
\end{aligned}
\]
The nonnegativity of the integrand therefore implies
\[
\int_0^{r/s-1} v^{\gamma-1}(1+v)^{-\gamma}\,dv \leq\log\frac
rs+\Theta(\gamma).
\]
This proves (\ref{eq:beta-claim}) and completes the proof.
\end{proof}

\begin{remark}\label{rem:radial-beta}
If $u$ is radially non-increasing, then for every $s\in J$ the set
$\{u>u^*(s)\}$ is the centered ball of measure $s$, up to a null
set. Consequently, $z(s)=\omega_n^{(1-\kappa)/n}s^{\gamma-1}$ and
$A(s)=A_0s^\gamma$ on $J$, while $\beta=0$ off $J$.

Thus
\[
\beta(s)=\frac{\chi_J(s)}{ns}\leq\frac1{ns},\qquad
\exp\left(c\int_s^r\beta(\tau)\,d\tau\right)\leq(r/s)^{c/n}.
\]
The radial version of (\ref{eq:pointwise}) therefore holds with
constant $1/(n\omega_n^{(1-\kappa)/n})$. The integrating factor is
exactly $(r/s)^{c/n}$ if $(u^*)'<0$ almost everywhere on $(0,S)$;
the absence of constant intervals alone does not imply this.

When $c>0$, (\ref{eq:Theta-bound}) gives
$e^{c\Theta(\gamma)/n}\leq e^{c/(n-1+\kappa)}$, with strict
inequality, and this factor tends to one as $\kappa\to1$.
\end{remark}

We shall combine the pointwise estimate of Lemma~\ref{lem:pointwise}
with a classical rearrangement inequality. For every nonnegative
measurable function $\varphi$ on $\R^n$ and $d\geq0$, the
Hardy--Littlewood inequality and the Hardy--Littlewood--P\'olya
principle (see, for example, \cite[Chapter~2]{BS}) yield
\begin{equation}\label{eq:HL-power}
\big\||x|^{-d}\varphi\big\|_Y \leq \omega_n^{d/n}
\big\|t^{-d/n}\varphi^*(t)\big\|_{\Ybar}.
\end{equation}
Indeed, for every measurable set $E$ with $|E|=t$,
\[
\int_E |x|^{-d}\varphi(x)\,dx \leq \omega_n^{d/n} \int_0^t
s^{-d/n}\varphi^*(s)\,ds,
\]
which gives the required submajorization. Equality holds in
(\ref{eq:HL-power}) when $\varphi$ is radially non-increasing.

We now use the pointwise estimate of Lemma~\ref{lem:pointwise},
together with the rearrangement inequality (\ref{eq:HL-power}), to
establish our main theorem.

\begin{theorem}
\label{thm:main} Let $n\geq2$, $b\geq0$, $a\leq b$, and let
$X,Y$ be r.i.\ Banach function spaces on $\R^n$, with representation
spaces $\Xbar,\Ybar$ on $(0,\infty)$. The following conditions are
equivalent:
\begin{align}
\big\||x|^{-b}f\big\|_Y &\leq C\big\||x|^{-a}\nabla f\big\|_X,
&&f\in C_c^\infty(\R^n);
\label{eq:CKN}\\
\|\HH_{a,b}h\|_{\Ybar} &\leq C_H\|h\|_{\Xbar}, &&0\leq h\in\Xbar.
\label{eq:bounded}
\end{align}
Here $\HH_{a,b}$ is defined in (\ref{eq:Hardy-operator}). If either
condition holds, then $t^{-b/n}\chi_{(0,1)}\in\Ybar$, or
equivalently $|x|^{-b}\chi_{\{|x|<1\}}\in Y$; in particular, $b<n$.
Moreover, $a\leq n-1$.
\end{theorem}

\begin{proof}
Assume first that (\ref{eq:bounded}) holds. We begin by deriving the
restrictions on $a$ and $b$. Taking $h=\chi_{(1,2)}$, we have, for
$0<t<1$,
\[
\HH_{a,b}h(t)=c_{a,n}t^{-b/n}, \qquad c_{a,n}:=\int_1^2
s^{(a+1)/n-1}\,ds>0.
\]
It follows that $t^{-b/n}\chi_{(0,1)}\in\Ybar$. By the local
integrability axiom for Banach function spaces,
\[
\int_0^1 t^{-b/n}\,dt<\infty,
\]
which forces $b<n$. Next, for $R>1$, applying (\ref{eq:bounded}) to
$\chi_{(R,R+1)}$ gives
\[
\begin{aligned}
R^{(a+1)/n-1} \big\|t^{-b/n}\chi_{(0,1)}\big\|_{\Ybar} &\preceq
\|\HH_{a,b}\chi_{(R,R+1)}\|_{\Ybar}\\
&\leq C_H\|\chi_{(R,R+1)}\|_{\Xbar}\\
&=C_H\|\chi_{(0,1)}\|_{\Xbar}.
\end{aligned}
\]
As the constants are independent of $R$, letting $R\to\infty$ yields
$a\leq n-1$.

We now establish (\ref{eq:CKN}). Fix $f\in C_c^\infty(\R^n)$ and put
$g=|x|^{-a}|\nabla f|$. We may assume $\|g\|_X<\infty$, since
otherwise there is nothing to prove. Choose $c$ such that
\begin{equation}\label{eq:c-range}
\max\{a,0\}\leq c\leq b, \qquad c-a\leq1\ \text{ if } c>0.
\end{equation}
This is always possible: if $a\leq0$ one may take $c=0$, and if
$a>0$ every $c\in[a,\min\{b,a+1\}]$ is admissible. When $b-a\leq1$,
(\ref{eq:c-range}) reduces to $\max\{a,0\}\leq c\leq b$. Set
\[
\kappa=c-a, \qquad u=|x|^{-c}|f|.
\]
Then $\kappa\geq0$, $b-c\geq0$, and $\kappa\leq1$ whenever $c>0$. By the product rule
and the inequality $|\nabla|f||\leq|\nabla f|$, we have
\begin{equation}\label{eq:grad-u}
|\nabla u| \leq |x|^{-c}|\nabla f|+c|x|^{-1}u \qquad\text{for a.e.
}x\ne0.
\end{equation}

For $M>0$, set $u_M=\min\{u,M\}$. We claim that
\begin{equation}\label{eq:claim-W11}
u_M\in W^{1,1}(\R^n).
\end{equation}
Assuming this for the moment, the truncation rule gives
\[
\nabla u_M=\chi_{\{u<M\}}\nabla u \qquad\text{a.e. in }\R^n.
\]
Since $u_M=u$ on $\{u<M\}$ and $|x|^\kappa|x|^{-c}=|x|^{-a}$,
inequality (\ref{eq:grad-u}) yields
\begin{equation*}
|x|^\kappa|\nabla u_M| \leq g+c|x|^{\kappa-1}u_M \qquad\text{a.e. in
}\R^n.
\end{equation*}
Choose a ball $B$ centered at the origin and containing the support
of $f$ in its interior. The local integrability property of Banach
function spaces ensures that $g\in L^1(B)$. We may therefore apply
Lemma~\ref{lem:pointwise} to $u_M$, with parameters $\kappa$ and
$c$, to obtain
\[
u_M^* \leq C_{n,\kappa,c}\, \HH_{a,c}(g_{*u_M}).
\]

Using the identity $t^{-(b-c)/n}\HH_{a,c}h(t)=\HH_{a,b}h(t)$, the
rearrangement estimate (\ref{eq:HL-power}) with $d=b-c$, the
boundedness of $\HH_{a,b}$, and the contraction property
(\ref{eq:rel-contract}), we obtain
\[
\begin{aligned}
\big\||x|^{-(b-c)}u_M\big\|_Y &\leq \omega_n^{(b-c)/n}
\big\|t^{-(b-c)/n}u_M^*(t)\big\|_{\Ybar}\\
&\leq \omega_n^{(b-c)/n}C_{n,\kappa,c}
\|\HH_{a,b}(g_{*u_M})\|_{\Ybar}\\
&\leq K_c\,C_H\|g_{*u_M}\|_{\Xbar}\\
&\leq K_c\,C_H\|g\|_X,
\end{aligned}
\]
where, by (\ref{eq:constant}) and with $\delta=b-a$,
\begin{equation}\label{eq:Kc}
K_c =\omega_n^{(b-c)/n}C_{n,c-a,c}
=\frac{E_c}{n\omega_n^{(1-\delta)/n}},
\end{equation}
where $\gamma_c=(n-1+c-a)/n$ and
\[
E_c=\begin{cases}
1, & c=0,\\
\exp(c\,\Theta(\gamma_c)/n), & c>0.
\end{cases}
\]
The bound is independent of $M$, while
$|x|^{-(b-c)}u_M\uparrow|x|^{-b}|f|$ almost everywhere as
$M\to\infty$. The Fatou property therefore gives
\begin{equation}\label{eq:CKN-c}
\big\||x|^{-b}f\big\|_Y \leq K_c\,C_H \big\||x|^{-a}\nabla
f\big\|_X.
\end{equation}

It remains to verify (\ref{eq:claim-W11}). The case $c=0$ is
immediate. For $c>0$, the function $u_M$ is locally Lipschitz away
from the origin, and (\ref{eq:grad-u}) gives
\[
|\nabla u_M| \leq \|\nabla f\|_\infty |x|^{-c} +cM|x|^{-1}
\qquad\text{a.e. on }\R^n\setminus\{0\}.
\]
As $c\leq b<n$ and $n\geq2$, the right-hand side is locally
integrable. Together with the compact support and boundedness of
$u_M$, this shows that both $u_M$ and its gradient away from the
origin belong to $L^1(\R^n)$. To identify the weak gradient, choose
$R>0$ such that $\supp f\Subset B(0,R)$ and integrate by parts on
$D_\varepsilon=B(0,R)\setminus\overline{B(0,\varepsilon)}$. For
$\varphi\in C_c^\infty(\R^n)$ and each coordinate $i$,
\[
\int_{D_\varepsilon}u_M\partial_i\varphi\,dx
=-\int_{D_\varepsilon}(\partial_i u_M)\varphi\,dx +\int_{\partial
B(0,\varepsilon)}u_M\varphi\nu_i\,dS,
\]
where $\nu$ denotes the outward normal to $D_\varepsilon$; the outer
boundary term vanishes. The remaining boundary integral is bounded
in absolute value by $n\omega_n
M\|\varphi\|_\infty\varepsilon^{n-1}$ and therefore tends to zero.
Since the volume integrals converge by integrability, we obtain
(\ref{eq:claim-W11}).

For the converse implication, suppose that (\ref{eq:CKN}) holds with
constant $C$. Given $0\leq h\in C_c^\infty(0,\infty)$, define
\[
F(r)=\int_r^\infty \rho^a h(\omega_n\rho^n)\,d\rho, \qquad
f(x)=F(|x|).
\]
Then $f\in C_c^\infty(\R^n)$ and is constant near the origin. Direct
differentiation and the change of variables $t=\omega_n r^n$ give
\[
|x|^{-a}|\nabla f(x)|=h(\omega_n|x|^n)
\]
and
\[
|x|^{-b}f(x) =\frac{\omega_n^{(b-a-1)/n}}{n}
\HH_{a,b}h(\omega_n|x|^n).
\]
Since the radial function $x\mapsto\psi(\omega_n|x|^n)$ is
equimeasurable with $\psi$ on $(0,\infty)$, inequality
(\ref{eq:CKN}) yields
\begin{equation}\label{eq:radial}
\|\HH_{a,b}h\|_{\Ybar} \leq n\omega_n^{(a+1-b)/n}C\|h\|_{\Xbar}.
\end{equation}

We extend (\ref{eq:radial}) to arbitrary nonnegative $h\in\Xbar$
without assuming density of smooth functions. For $j\geq2$, set
$h_j=\min\{h,j\}\chi_{(1/j,j)}$. Then $h_j\uparrow h$, $h_j\in
L^1(0,\infty)$, and $\|h_j\|_{\Xbar}\leq\|h\|_{\Xbar}$. Extend $h_j$
by zero to $\R$ and consider
$h_{j,\varepsilon}=h_j*\rho_\varepsilon$, where
$\rho_\varepsilon\geq0$ is a smooth approximate identity supported
in $(-\varepsilon,\varepsilon)$. For $\varepsilon<1/(2j)$, we have
$h_{j,\varepsilon}\in C_c^\infty(0,\infty)$, with $\supp
h_{j,\varepsilon}\subset K_j=[1/(2j),j+1]$.

Mollification is contractive by submajorization. Indeed, for every
measurable set $E\subset\R$ of finite measure, Fubini's theorem and
the Hardy--Littlewood inequality give
\[
\begin{aligned}
\int_E h_{j,\varepsilon}(x)\,dx &=\int_\R \rho_\varepsilon(y)
  \int_{E-y}h_j(x)\,dx\,dy\\
&\leq\int_0^{|E|}h_j^*(s)\,ds.
\end{aligned}
\]
Taking the supremum over sets of fixed measure, we obtain
$h_{j,\varepsilon}\prec\!\prec h_j$, and the
Hardy--Littlewood--P\'olya principle (\ref{eq:HLP}) therefore
implies
\[
\|h_{j,\varepsilon}\|_{\Xbar} \leq\|h_j\|_{\Xbar} \leq\|h\|_{\Xbar}.
\]
On the other hand, $h_{j,\varepsilon}\to h_j$ in $L^1(\R)$, and the
kernel $s^{(a+1)/n-1}$ is bounded on $K_j$. Thus, for every $t>0$,
\[
\begin{aligned}
&|\HH_{a,b}h_{j,\varepsilon}(t)
  -\HH_{a,b}h_j(t)|\\
&\qquad\leq t^{-b/n}\sup_{s\in K_j}s^{(a+1)/n-1}
\|h_{j,\varepsilon}-h_j\|_1 \longrightarrow0.
\end{aligned}
\]
Applying (\ref{eq:radial}) to $h_{j,\varepsilon}$ and using the
Fatou property as $\varepsilon\downarrow0$, we conclude that
\[
\|\HH_{a,b}h_j\|_{\Ybar} \leq n\omega_n^{(a+1-b)/n} C\|h\|_{\Xbar}.
\]
Finally, the positivity of $\HH_{a,b}$ ensures that
$\HH_{a,b}h_j\uparrow\HH_{a,b}h$. A second application of the Fatou
property proves (\ref{eq:bounded}).

This completes the proof.
\end{proof}

The restriction $a\leq b$ in Theorem~\ref{thm:main} is necessary,
for arbitrary real exponents.

\begin{proposition}\label{prop:a-le-b}
Let $n\geq1$, $a,b\in\R$, and let $X,Y$ be r.i.\ Banach function
spaces on $\R^n$. If (\ref{eq:CKN}) holds, then $a\leq b$.
\end{proposition}

\begin{proof}
Fix $0\neq\varphi\in C_c^\infty(B(0,1))$ and, for $R\geq2$, put
$\varphi_R(x)=\varphi(x-Re_1)$. Since $R/2\leq|x|\leq 2R$ on
$\supp\varphi_R$, the lattice property and the translation
invariance of r.i.\ norms give
\[
2^{-|b|}R^{-b}\|\varphi\|_Y \leq\big\||x|^{-b}\varphi_R\big\|_Y \leq
C\big\||x|^{-a}\nabla\varphi_R\big\|_X \leq
C\,2^{|a|}R^{-a}\big\|\nabla\varphi\big\|_X.
\]
Since $\|\varphi\|_Y>0$ and $\|\nabla\varphi\|_X<\infty$, letting
$R\to\infty$ yields $a\leq b$.
\end{proof}

\begin{remark}\label{rem:effective-range}
The equivalence in Theorem~\ref{thm:main} holds for all $b\geq0$ and
$a\leq b$, but the existence of an admissible pair of r.i.\ spaces
$X,Y$ requires, in addition,
\[
\delta=b-a<n+1.
\]
Indeed, since necessarily $b<n$, if $\delta\geq n+1$, then $a+1<0$.
For $h=\chi_{(0,1)}$ and $0<t<1$,
\[
\HH_{a,b}h(t) =\frac{t^{(1-\delta)/n}-t^{-b/n}}{-(a+1)/n} \simeq
t^{(1-\delta)/n} \qquad (t\downarrow0).
\]
This function is not locally integrable when $\delta\geq n+1$,
contradicting the boundedness of $\HH_{a,b}$ between Banach function
spaces.

Conversely, if $0\leq b<n$, $a\leq\min\{b,n-1\}$ and $\delta<n+1$,
admissible pairs exist. Indeed, for $X=L^1\cap L^\infty$ we have
$\Xbar'=L^1+L^\infty$, which contains every function comparable to
$s^{(1-\delta)/n}$ on $(0,1)$ and to $s^{(a+1)/n-1}$ on
$(1,\infty)$; hence Theorem~\ref{thm:optimal}(i) provides an r.i.\
target.
\end{remark}

\subsection{A perturbative argument and its limitations}
\label{subsec:perturbative}

When $a>0$ and a weighted Hardy inequality is available in $X$, the
sufficiency part of Theorem~\ref{thm:main} can be obtained without
the integrating-factor argument. We record this fact, since it
delimits precisely where the method of Lemma~\ref{lem:pointwise} is
needed.

\begin{proposition}\label{prop:perturbative}
Let $n\geq2$, $0<a\leq b$, $\delta=b-a$, and let $X,Y$ be r.i.\
Banach function spaces on $\R^n$. Assume that:
\begin{itemize}
\item[(i)] there is $C_1>0$ such that
\begin{equation}\label{eq:hardy-X}
\big\||x|^{-a-1}f\big\|_X\leq C_1\big\||x|^{-a}\nabla f\big\|_X,
\qquad f\in C_c^\infty(\R^n);
\end{equation}
\item[(ii)] $\HH_{0,\delta}:\Xbar\to\Ybar$ is bounded, with norm
$C_2$.
\end{itemize}
Then (\ref{eq:CKN}) holds, with
$C\leq(1+aC_1)C_2/(n\omega_n^{(1-\delta)/n})$. Moreover,
$\HH_{0,\delta}h\leq\HH_{a,b}h$ for every $h\geq0$, so that
\textup{(ii)} holds whenever (\ref{eq:bounded}) does.
\end{proposition}

\begin{proof}
By Theorem~\ref{thm:main}, applied with $b$ replaced by $a+1$,
hypothesis (i) implies $a+1<n$. Fix $f\in C_c^\infty(\R^n)$, put
$g=|x|^{-a}|\nabla f|$, and assume $\|g\|_X<\infty$. Let
$u=|x|^{-a}|f|$ and $u_M=\min\{u,M\}$. As in the proof of
Theorem~\ref{thm:main} (with $c=a$), $u_M$ is a compactly supported
function in $W^{1,1}(\R^n)$ and, by (\ref{eq:grad-u}),
\[
|\nabla u_M|\leq\chi_{\{u<M\}}\bigl(g+a|x|^{-a-1}|f|\bigr)
\quad\text{a.e. in }\R^n.
\]
Hence, by (i), $\||\nabla u_M|\|_X\leq(1+aC_1)\|g\|_X$. We apply
Lemma~\ref{lem:pointwise} to $u_M$ with $c=0$, $\kappa=0$ and with
$|\nabla u_M|$ in place of $g$; in this case no zeroth-order term is
present and the lemma reduces to the classical pointwise estimate
\[
u_M^*\leq\frac1{n\omega_n^{1/n}}\,
\HH_{0,0}\bigl((|\nabla u_M|)_{*u_M}\bigr).
\]
Using (\ref{eq:HL-power}) with $d=\delta$, the identity
$t^{-\delta/n}\HH_{0,0}=\HH_{0,\delta}$, hypothesis (ii) and
(\ref{eq:rel-contract}), we obtain
\[
\big\||x|^{-\delta}u_M\big\|_Y
\leq\frac{C_2}{n\omega_n^{(1-\delta)/n}}\big\||\nabla u_M|\big\|_X
\leq\frac{(1+aC_1)C_2}{n\omega_n^{(1-\delta)/n}}\,\|g\|_X .
\]
Since $|x|^{-\delta}u_M\uparrow|x|^{-b}|f|$ as $M\to\infty$, the
Fatou property gives (\ref{eq:CKN}). Finally, for $h\geq0$ and
$a\geq0$,
\[
\HH_{a,b}h(t)=t^{-\delta/n}\int_t^\infty h(s)
\Bigl(\frac st\Bigr)^{a/n}s^{1/n-1}\,ds\geq\HH_{0,\delta}h(t).
\qedhere
\]
\end{proof}

\begin{remark}\label{rem:perturbative}
(a) By Theorem~\ref{thm:main}, (\ref{eq:hardy-X}) is equivalent to
the boundedness of $\HH_{a,a+1}$ on $\Xbar$. For $X=L^p$ this holds
if and only if $a<n/p-1$; hence
Proposition~\ref{prop:perturbative} recovers the sufficiency part of
the subcritical Lebesgue inequalities
(Proposition~\ref{cor:lebesgue}) using only rearrangements with
respect to Lebesgue measure and the unweighted isoperimetric
inequality. A sufficient condition for (\ref{eq:hardy-X}) in terms
of the Boyd index of $X$ is given in Proposition~\ref{cor:hardy-ri}.

(b) Hypothesis (i) is neither necessary for (\ref{eq:CKN}) nor
implied by (\ref{eq:bounded}). Let $1\leq p<n$, $a=n/p-1>0$, $a\leq
b<n$ and $1/q=b/n$. By Proposition~\ref{cor:critical},
$\HH_{a,b}:L^{p,1}\to L^{q,\infty}$ is bounded and (\ref{eq:CKN})
holds for $X=L^{p,1}$, $Y=L^{q,\infty}$. On the other hand,
(\ref{eq:hardy-X}) fails for $X=L^{p,1}$: by
Theorem~\ref{thm:main} it would force
$t^{-(a+1)/n}\chi_{(0,1)}=t^{-1/p}\chi_{(0,1)}\in L^{p,1}(0,\infty)$,
which is false. The same phenomenon occurs for every r.i.\ space $X$
such that $\HH_{a,b}$ is bounded from $\Xbar$ into some $\Ybar$
while $\HH_{a,a+1}$ is unbounded on $\Xbar$.

(c) Consequently, a criterion valid for all pairs $X,Y$ cannot be
obtained by absorbing the zeroth-order term in
(\ref{eq:weight-transfer}) through a norm inequality. In the proof
of Theorem~\ref{thm:main} this term is instead integrated exactly at
the level of rearrangements, by means of
Lemmas~\ref{lem:active} and~\ref{lem:pointwise}, and no property of
$X$ is used beyond the contraction (\ref{eq:rel-contract}).
\end{remark}

\subsection{Optimal constants}
\label{subsec:optimal-constants}

\begin{proposition}
\label{prop:sharp} Under the assumptions of Theorem~\ref{thm:main},
suppose that (\ref{eq:CKN}) holds, and put $\delta=b-a$. Let
$C_{\mathrm{CKN}}^{\mathrm{opt}}$ denote the optimal constant in
(\ref{eq:CKN}), and let $C_{\mathrm{rad}}^{\mathrm{opt}}$ denote the
optimal constant when (\ref{eq:CKN}) is restricted to radial
functions. Then
\begin{equation}
\label{eq:radial-sharp} C_{\mathrm{rad}}^{\mathrm{opt}} =
\frac{1}{n\omega_n^{(1-\delta)/n}} \|\HH_{a,b}\|_{\Xbar\to\Ybar}.
\end{equation}

Moreover,
\begin{equation}
\label{eq:radial-comparison} C_{\mathrm{rad}}^{\mathrm{opt}} \leq
C_{\mathrm{CKN}}^{\mathrm{opt}} \leq e^{\Lambda(n,a,b)}
C_{\mathrm{rad}}^{\mathrm{opt}},
\end{equation}
where
\begin{equation}\label{eq:Lambda}
\Lambda(n,a,b)=
\begin{cases}
0,&a\leq0\ \text{or}\ \delta\geq1,\\[3pt]
\displaystyle\frac1n\min\left\{ a\Theta\left(\frac{n-1}{n}\right),\,
b\Theta\left(\frac{n-1+\delta}{n}\right)\right\},&a>0,\ \delta<1.
\end{cases}
\end{equation}
In all cases $e^{\Lambda(n,a,b)}\leq2$, and therefore
\begin{equation}\label{eq:uniform-comparison}
C_{\mathrm{rad}}^{\mathrm{opt}} \leq C_{\mathrm{CKN}}^{\mathrm{opt}}
\leq2C_{\mathrm{rad}}^{\mathrm{opt}}.
\end{equation}

In particular, if $a\leq0$ or $\delta\geq1$, then
\begin{equation}
\label{eq:sharp} C_{\mathrm{CKN}}^{\mathrm{opt}} =
C_{\mathrm{rad}}^{\mathrm{opt}} =
\frac{1}{n\omega_n^{(1-\delta)/n}}\|\HH_{a,b}\|_{\Xbar\to\Ybar}.
\end{equation}
\end{proposition}

\begin{proof}

The necessity argument in the proof of Theorem~\ref{thm:main} uses
only radial test functions. Hence
\[
\|\HH_{a,b}\|_{\Xbar\to\Ybar} \leq n\omega_n^{(1-\delta)/n}
C_{\mathrm{rad}}^{\mathrm{opt}}.
\]

For the reverse inequality, let $f\in C_c^\infty(\R^n)$ be radial,
with $f(x)=F(|x|)$. Since $F$ vanishes for sufficiently large $r$,
the fundamental theorem of calculus gives
\[
|F(r)|\leq\int_r^\infty |F'(\rho)|\,d\rho.
\]

Define $h$ on $(0,\infty)$ by
\[
h(\omega_n r^n)=r^{-a}|F'(r)|.
\]
A direct change of variables gives
\[
r^{-b}\int_r^\infty |F'(\rho)|\,d\rho =
\frac{1}{n\omega_n^{(1-\delta)/n}} \HH_{a,b}h(\omega_n r^n).
\]

Moreover, $x\mapsto h(\omega_n|x|^n)$ is equimeasurable with $h$ and
\[
h(\omega_n|x|^n)=|x|^{-a}|\nabla f(x)|.
\]
Therefore
\[
\begin{aligned}
\big\||x|^{-b}f\big\|_Y &\leq \frac{1}{n\omega_n^{(1-\delta)/n}}
\|\HH_{a,b}h\|_{\Ybar}\\
&\leq \frac{1}{n\omega_n^{(1-\delta)/n}}
\|\HH_{a,b}\|_{\Xbar\to\Ybar} \big\||x|^{-a}\nabla f\big\|_X.
\end{aligned}
\]
This proves (\ref{eq:radial-sharp}).

Since radial functions form a subclass of $C_c^\infty(\R^n)$,
\[
C_{\mathrm{rad}}^{\mathrm{opt}} \leq
C_{\mathrm{CKN}}^{\mathrm{opt}}.
\]

On the other hand, by (\ref{eq:CKN-c}), (\ref{eq:Kc}) and
(\ref{eq:radial-sharp}),
\[
C_{\mathrm{CKN}}^{\mathrm{opt}} \leq
K_c\|\HH_{a,b}\|_{\Xbar\to\Ybar}
=E_c\,C_{\mathrm{rad}}^{\mathrm{opt}}
\]
for every $c$ satisfying (\ref{eq:c-range}).

If $a\leq0$, the choice $c=0$ gives equality. If $a>0$ and
$\delta\geq1$, the choice $c=a+1$ is admissible in
(\ref{eq:c-range}) and gives $\gamma_c=1$, hence
$\Theta(\gamma_c)=0$. This proves (\ref{eq:sharp}).

Suppose now that $a>0$ and $\delta<1$, so that the admissible values
of $c$ form the interval $[a,b]$. By Theorem~\ref{thm:main}, we have $a\leq
n-1$, so that
\[
d=\frac{n-1-a}{n}\geq0.
\]
Consider the function $G(t)=t\Theta(d+t)$, defined for $0<t<1-d$.
Then $E_c=\exp(G(c/n))$, so it suffices to minimize $G$ on
$[a/n,b/n]$. Termwise differentiation of the series in (\ref{eq:Theta}) gives
\[
G''(t)=-2\sum_{k=0}^\infty \frac{k+d}{(k+d+t)^3}\leq0.
\]

Thus $G$ is concave, and its minimum on $[a/n,b/n]$ is attained at
one of the endpoints. These correspond to $c=a$ and $c=b$,
respectively, and yield (\ref{eq:Lambda}) and
(\ref{eq:radial-comparison}).

It remains to establish the uniform bound. We may assume $a>0$ and
$\delta<1$, since otherwise $\Lambda=0$. Set $\gamma_0=(n-1)/n\in[1/2,1)$. As
$a\leq n-1$,
\[
\Lambda \leq\frac an\Theta(\gamma_0) \leq\gamma_0\Theta(\gamma_0).
\]
For $1/2\leq\gamma\leq1$, differentiation of the series in
(\ref{eq:Theta}) gives
\[
\frac{d}{d\gamma} \bigl[\gamma\Theta(\gamma)\bigr]
=\sum_{k=0}^\infty \left( \frac{k}{(k+\gamma)^2}-\frac1{k+1}
\right)\leq0,
\]
since
\[
(k+\gamma)^2-k(k+1) =(2\gamma-1)k+\gamma^2\geq0.
\]
Consequently,
\[
\gamma_0\Theta(\gamma_0) \leq\frac12\Theta\left(\frac12\right).
\]
By (\ref{eq:Theta}),
\[
\Theta\left(\frac12\right) =\sum_{k=0}^\infty
\left(\frac{2}{2k+1}-\frac{2}{2k+2}\right) =2\log2.
\]
Therefore $\Lambda\leq\log2$, so that $e^{\Lambda(n,a,b)}\leq2$,
proving (\ref{eq:uniform-comparison}).
\end{proof}

\begin{remark}\label{rem:Lambda}\label{rem:symmetry}
The choices $c=a$ and $c=b$ in (\ref{eq:Lambda}) correspond,
respectively, to transferring the gradient weight or the entire
function weight to $u$. For $a>0$ and $\delta<1$, we have
\[
\Lambda \leq \frac{b}{n} \Theta\left(\frac{n-1+\delta}{n}\right)
\leq \frac{b}{n-1+\delta}.
\]
In dimension two, the alternative choice $c=a$ gives $e^\Lambda\leq
2^a$. The universal factor $2$ in (\ref{eq:uniform-comparison}) is
therefore only an upper bound, and we make no claim concerning its
optimality.

For the Lebesgue inequalities of Proposition~\ref{cor:lebesgue} with
$p>1$, the general and radial optimal constants coincide throughout
the present range. This follows from \cite[Section~8.2, particularly
Proposition~8.1 and (8.22)]{AlvinoEtAl}, after reversing the signs
of the weight exponents; see also \cite{ChouChu,Horiuchi} for $a\geq0$. Indeed, their condition $a\leq a_2$
corresponds to our assumption $b\geq0$. For symmetry and symmetry
breaking in other ranges of the classical CKN inequalities, see, for
example, \cite{CatrinaWang,FelliSchneider,DolbeaultEstebanLoss}.

Thus the factor $e^\Lambda$ reflects a loss in the pointwise
argument, at least in the Lebesgue case. Whether the general and
radial optimal constants coincide for arbitrary r.i.\ pairs when
$a>0$ and $\delta<1$ is not settled by the present argument.
\end{remark}

\begin{remark}\label{rem:tests}\label{rem:radial}
The test class in (\ref{eq:CKN}) includes functions that are
constant near the origin. Such tests force
$|x|^{-b}\chi_{\{|x|<1\}}\in Y$ and distinguish our setting from
inequalities posed only on $C_c^\infty(\R^n\setminus\{0\})$. The
case $a=b=0$ recovers the unweighted reduction, with
\[
\HH_{0,0}h(t) =\int_t^\infty h(s)s^{1/n-1}\,ds.
\]
When $b<0$, the weight $|x|^{-b}$ is radially increasing, and the
rearrangement estimate (\ref{eq:HL-power}) no longer applies in the
required form.

For radial functions, however, neither $b\geq0$ nor $a\leq b$ is
needed. The radial arguments in the proofs of
Theorem~\ref{thm:main} and Proposition~\ref{prop:sharp} apply to
arbitrary real $a,b$, yielding the same equivalence with
Hardy-operator boundedness and the optimal constant
(\ref{eq:radial-sharp}). For general functions, $a\leq b$ is
necessary by Proposition~\ref{prop:a-le-b}.

In Lebesgue spaces, the corresponding weighted Hardy criterion
applies whenever $1\leq p\leq q<\infty$ and
\[
\frac1q=\frac1p-\frac{1-\delta}{n},
\]
giving boundedness precisely when $1/q>b/n$. These conditions allow
$-(n-1)<\delta\leq1$, including negative values of $\delta$. For
related radial Caffarelli--Kohn--Nirenberg inequalities in the range
$\delta<0$, see \cite[Theorem~18]{MallickNguyen}.
\end{remark}

\section{Optimal target spaces}
\label{sec:optimal}

Theorem~\ref{thm:main} reduces the search for the optimal r.i.\
target space in (\ref{eq:CKN}) to the corresponding problem for
$\HH_{a,b}$. Since $b\geq0$, both $\HH_{a,b}h$, for $h\geq0$, and
the weight $t^{-b/n}$ are non-increasing. This allows us to
construct the optimal target by a direct argument based on associate
spaces. For general weighted Hardy operators, optimal domains and
targets are studied in \cite{Mihula}; for unweighted Sobolev
embeddings on $\R^n$, see, for example, \cite{ACPS,MihulaPRSE}.

We use the associate-space duality introduced in
Section~\ref{sec:preliminaries}. Set
\[
\eta=\frac{a+1}{n}-1, \qquad \HH'_{a,b}g(s) =s^\eta\int_0^s
g(t)t^{-b/n}\,dt.
\]
The operator $\HH'_{a,b}$ is the formal adjoint of $\HH_{a,b}$ with
respect to the integral pairing. Indeed, Tonelli's theorem gives
\begin{equation}\label{eq:adjoint}
\int_0^\infty \HH_{a,b}h(t)g(t)\,dt =\int_0^\infty
h(s)\HH'_{a,b}g(s)\,ds
\end{equation}
for all nonnegative measurable functions $h,g$.

In the notation of \cite[Proposition~3.7]{Mihula}, the operator
$\HH_{a,b}$ corresponds to
\[
u(t)=t^{-b/n}, \qquad v(t)=t^\eta, \qquad \nu(t)=t.
\]
The auxiliary function $\xi$ defined in \cite[(3.1)]{Mihula} takes
the form
\[
\xi(s)= \frac{n}{n-b}s^{(1-\delta)/n}\chi_{(0,1)}(s)
+s^\eta\chi_{[1,\infty)}(s).
\]
If $0\leq\delta\leq1$, the first term belongs to $\Xbar'$, and
therefore $\xi\in\Xbar'$ is equivalent to
\[
s^\eta\chi_{(1,\infty)}(s)\in\Xbar'.
\]
If $\delta>1$, the first term is unbounded near the origin, and its
membership in $\Xbar'$ is a genuine condition. In all cases, since
$\HH'_{a,b}(\chi_{(0,1)})(s)=\frac{n}{n-b}s^\eta\min\{s,1\}^{1-b/n}$
is comparable to $\xi$, the condition $\xi\in\Xbar'$ can also be
expressed as $\HH'_{a,b}(\chi_{(0,1)})\in\Xbar'$. Thus the optimal
target construction follows from the general result cited above. We
include a direct proof in our setting to make explicit its
application to the CKN inequality through Theorem~\ref{thm:main}.

\begin{theorem}\label{thm:optimal}
Let $n\geq2$, $0\leq b<n$, $\delta=b-a\geq0$, and let $X$ be an
r.i.\ space on $\R^n$.
\begin{itemize}
\item[(i)] Inequality (\ref{eq:CKN}) holds for some r.i.\ space $Y$
if and only if $\HH'_{a,b}(\chi_{(0,1)})\in\Xbar'$. If
$\delta\leq1$, this is equivalent to
$s^\eta\chi_{(1,\infty)}(s)\in\Xbar'$.
\item[(ii)] Assume that $\HH'_{a,b}(\chi_{(0,1)})\in\Xbar'$. Then
\begin{equation*}
\rho(g)=\big\|\HH'_{a,b}(g^*)\big\|_{\Xbar'}
\end{equation*}
is an r.i.\ Banach function norm on $(0,\infty)$. Let $\Ybar_X$ be
the associate space of the space determined by $\rho$, and let $Y_X$
be the r.i.\ space on $\R^n$ with representation space $\Ybar_X$.
Then (\ref{eq:CKN}) holds with $Y=Y_X$, and $Y_X\hookrightarrow Y$
for every r.i.\ space $Y$ for which (\ref{eq:CKN}) holds.
\end{itemize}
\end{theorem}

\begin{proof}
(i) Suppose that (\ref{eq:CKN}) holds for some $Y$. By
Theorem~\ref{thm:main}, $\HH_{a,b}:\Xbar\to\Ybar$ is bounded, with
norm $C_H$. By the local integrability axiom, there is $C_1$ such
that $\int_0^1|k|\leq C_1\|k\|_{\Ybar}$ for every $k\in\Ybar$. For
$0\leq h\in\Xbar$, (\ref{eq:adjoint}) gives
\[
\int_0^\infty h(s)\HH'_{a,b}(\chi_{(0,1)})(s)\,ds
=\int_0^1\HH_{a,b}h(t)\,dt\leq C_1C_H\|h\|_{\Xbar}.
\]
Thus $\HH'_{a,b}(\chi_{(0,1)})\in\Xbar'$. The converse follows from
(ii).

(ii) We verify that $\rho$ is an r.i.\ Banach function norm in the
sense of \cite[Chapter~1]{BS}. Since $\rho(g)$ depends only on
$g^*$, it is rearrangement-invariant. Positive homogeneity follows
from the linearity of $\HH'_{a,b}$, while $|g_1|\leq|g_2|$ implies
$g_1^*\leq g_2^*$ and hence $\rho(g_1)\leq\rho(g_2)$.

\emph{Subadditivity.} For $s>0$, the function
$t^{-b/n}\chi_{(0,s)}(t)$ is non-increasing, since $b\geq0$. As
$(g_1+g_2)^*\prec\!\prec g_1^*+g_2^*$ \cite[Chapter~2]{BS}, Hardy's
lemma gives
\[
\int_0^s(g_1+g_2)^*(t)t^{-b/n}\,dt \leq
\int_0^s\bigl(g_1^*(t)+g_2^*(t)\bigr)t^{-b/n}\,dt.
\]
Consequently,
\[
\HH'_{a,b}((g_1+g_2)^*) \leq \HH'_{a,b}(g_1^*)+\HH'_{a,b}(g_2^*),
\]
and the triangle inequality in $\Xbar'$ yields
\[
\rho(g_1+g_2)\leq\rho(g_1)+\rho(g_2).
\]

\emph{Fatou property.} If $0\leq g_k\uparrow g$ almost everywhere,
then $g_k^*\uparrow g^*$ almost everywhere. By monotone convergence,
\[
\HH'_{a,b}(g_k^*)\uparrow\HH'_{a,b}(g^*),
\]
and the Fatou property of $\Xbar'$ gives $\rho(g_k)\uparrow\rho(g)$.

\emph{Sets of finite measure.} For $r>0$,
\[
\HH'_{a,b}(\chi_{(0,r)})(s)
=\frac{n}{n-b}s^\eta\min\{s,r\}^{1-b/n}.
\]
Since $\min\{s,r\}\leq\max\{1,r\}\min\{s,1\}$ and $1-b/n>0$,
\[
\HH'_{a,b}(\chi_{(0,r)}) \leq\max\{1,r^{1-b/n}\}\,
\HH'_{a,b}(\chi_{(0,1)}),
\]
and the right-hand side belongs to $\Xbar'$ by assumption. Thus
\[
\rho(\chi_E)<\infty
\]
for every measurable set $E$ of finite measure.

\emph{Local integrability and definiteness.} Since $t^{-b/n}\geq
s^{-b/n}$ for $0<t<s$,
\[
\HH'_{a,b}(g^*)(s) \geq s^{\eta+1-b/n}g^{**}(s)
=s^{(1-\delta)/n}g^{**}(s).
\]
As $g^{**}$ is non-increasing, for every $r>0$, with
$m_r=\min_{r\leq s\leq2r}s^{(1-\delta)/n}>0$,
\[
\begin{aligned}
\rho(g) &\geq m_r\,g^{**}(2r)
\big\|\chi_{(r,2r)}\big\|_{\Xbar'}\\
&\geq \frac{m_r}{2r} \big\|\chi_{(r,2r)}\big\|_{\Xbar'} \int_0^r
g^*(t)\,dt.
\end{aligned}
\]
Hence, for every measurable set $E$ of measure $r$, the
Hardy--Littlewood inequality gives
\[
\begin{aligned}
\int_E|g(t)|\,dt
&\leq\int_0^r g^*(t)\,dt\\
&\leq \frac{2r} {m_r\|\chi_{(r,2r)}\|_{\Xbar'}}\,\rho(g).
\end{aligned}
\]
The constant is finite and positive. This proves the local
integrability axiom and also shows that $\rho(g)=0$ implies $g=0$
almost everywhere.

Thus $\rho$ satisfies all the axioms of an r.i.\ Banach function
norm. Denote by $\Zbar$ the corresponding space. By the
Lorentz--Luxemburg theorem, $\Ybar_X=\Zbar'$ is an r.i.\ Banach
function space, $\Ybar_X'=\Zbar$, and
\[
\|k\|_{\Ybar_X} = \sup\Bigl\{ \int_0^\infty
k^*(t)g^*(t)\,dt:\rho(g)\leq1 \Bigr\};
\]
see \cite[Chapters~1 and 2]{BS}.

\emph{Boundedness.} Let $0\leq h\in\Xbar$. Since $b\geq0$,
$\HH_{a,b}h$ is non-increasing, so $(\HH_{a,b}h)^*=\HH_{a,b}h$. By
(\ref{eq:adjoint}),
\[
\begin{aligned}
\|\HH_{a,b}h\|_{\Ybar_X} &=\sup_{\rho(g)\leq1}
\int_0^\infty\HH_{a,b}h(t)g^*(t)\,dt\\
&=\sup_{\rho(g)\leq1}
\int_0^\infty h(t)\HH'_{a,b}(g^*)(t)\,dt\\
&\leq\|h\|_{\Xbar}.
\end{aligned}
\]
Hence $\HH_{a,b}:\Xbar\to\Ybar_X$ is bounded, with norm at most one,
and (\ref{eq:CKN}) holds with $Y=Y_X$ by Theorem~\ref{thm:main}.

\emph{Optimality.} Suppose that (\ref{eq:CKN}) holds for some r.i.\
space $Y$. By Theorem~\ref{thm:main}, $\HH_{a,b}:\Xbar\to\Ybar$ is
bounded, with some norm $C$. For $g\geq0$, the associate norm
formula and (\ref{eq:adjoint}) give
\[
\begin{aligned}
\rho(g) &=\sup_{\substack{h\geq0\\\|h\|_{\Xbar}\leq1}}
\int_0^\infty h(t)\HH'_{a,b}(g^*)(t)\,dt\\
&=\sup_{\substack{h\geq0\\\|h\|_{\Xbar}\leq1}}
\int_0^\infty\HH_{a,b}h(t)g^*(t)\,dt\\
&\leq C\|g\|_{\Ybar'}.
\end{aligned}
\]
Thus $\Ybar'\hookrightarrow\Zbar=\Ybar_X'$, and the
Lorentz--Luxemburg theorem gives $\Ybar_X\hookrightarrow\Ybar$, that
is, $Y_X\hookrightarrow Y$.
\end{proof}

\section{Examples and applications}
\label{sec:examples}

Theorem~\ref{thm:main} transfers boundedness criteria and
optimal-space results for weighted Hardy operators to CKN
inequalities in arbitrary r.i.\ spaces; see, in particular,
\cite{Mihula}. This includes Orlicz and Lorentz--Zygmund spaces.
Here we focus on Lebesgue and Lorentz spaces, for which there is an
extensive literature
\cite{CKN,Yacoub,NursultanovTikhonov,LiuYuZhou,Wang}, and on the critical
endpoints and optimality properties supplied by the reduction
principle.

\subsection{Subcritical Lebesgue and Lorentz embeddings}

We first recover the classical Lebesgue range, including $p=1$. The
corresponding power-weight Hardy criterion is well known; see
\cite{CKN,KufnerMaligrandaPersson}.

\begin{proposition}\label{cor:lebesgue}
Let $n\geq2$, $b\geq0$, $0\leq\delta=b-a<1$, and $1\leq
p<n/(1-\delta)$. Define $q$ by $1/q=1/p-(1-\delta)/n$. Then
\begin{equation}\label{eq:lp}
\big\||x|^{-b}f\big\|_{L^q(\R^n)} \leq C\big\||x|^{-a}\nabla
f\big\|_{L^p(\R^n)}, \qquad f\in C_c^\infty(\R^n),
\end{equation}
holds if and only if $a<n/p-1$. In particular, $\delta=0$ recovers
the equal-weight case $a=b\geq0$.
\end{proposition}

\begin{proof}
Put $\lambda=1/p-(a+1)/n=1/q-b/n$. If (\ref{eq:lp}) holds,
Theorem~\ref{thm:main} gives the boundedness of $\HH_{a,b}:L^p\to
L^q$. Testing it on $h=\chi_{(1,2)}$ shows that
$t^{-b/n}\chi_{(0,1)}\in L^q$, so $bq<n$ and $\lambda>0$.

Conversely, suppose that $\lambda>0$. For $h\geq0$ set
$G(y)=e^{y/p}h(e^y)$. The changes of variables $t=e^x$ and
$s=e^{x+v}$ yield
\[
e^{x/q}\HH_{a,b}h(e^x) =\int_0^\infty e^{-\lambda v}G(x+v)\,dv.
\]
The changes of variables preserve the respective $L^p$ and $L^q$
norms. With $r=n/(n-1+\delta)>1$, we have $1+1/q=1/p+1/r$, and
Young's inequality gives
\[
\|\HH_{a,b}h\|_{L^q(0,\infty)} \leq
(r\lambda)^{-1/r}\|h\|_{L^p(0,\infty)}.
\]
Theorem~\ref{thm:main} completes the proof.
\end{proof}

The reduction also gives Lorentz refinements. In part~(i) below, the
target $L^{q,1}$ is equipped with the $f^*$-based norm defined in
(\ref{eq:Lp1}). In part~(ii), where no sharp constant is asserted,
either of the equivalent Lorentz norm conventions may be used.

\begin{corollary}\label{cor:lorentz-refined}
Let $n\geq2$, $b\geq0$, $\delta=b-a\geq0$, $1\leq p<\infty$,
$a<n/p-1$, and $1/q=1/p-(1-\delta)/n$.
\begin{itemize}
\item[(i)] If $p=1$ and $\delta<1$, then
\begin{equation*}
\big\||x|^{-b}f\big\|_{L^{q,1}(\R^n)} \leq\frac{e^{\Lambda(n,a,b)}}
{(n-1-a)\omega_n^{(1-\delta)/n}} \big\||x|^{-a}\nabla
f\big\|_{L^1(\R^n)}, \qquad f\in C_c^\infty(\R^n).
\end{equation*}
The radial optimal constant is $1/((n-1-a)\omega_n^{(1-\delta)/n})$;
it is also the general optimal constant when $a\leq0$.
\item[(ii)] If $p>1$ and $q>1$, then, for every
$1\leq r\leq\infty$,
\begin{equation*}
\big\||x|^{-b}f\big\|_{L^{q,r}(\R^n)} \preceq\big\||x|^{-a}\nabla
f\big\|_{L^{p,r}(\R^n)}, \qquad f\in C_c^\infty(\R^n).
\end{equation*}
In particular, $L^{q,p}$ is an admissible target for $X=L^p$. For
$\delta<1$ this refines the Lebesgue estimate; for $\delta>1$ it
replaces the inadmissible target $L^q$ (see
Remark~\ref{rem:delta-large}).
\end{itemize}
\end{corollary}

\begin{proof}
(i) Here $\lambda=1/q-b/n=(n-1-a)/n>0$. Since $\HH_{a,b}h$ is
non-increasing for $h\geq0$, Tonelli's theorem gives
\[
\|\HH_{a,b}h\|_{L^{q,1}} =\int_0^\infty t^{\lambda-1} \int_t^\infty
h(s)s^{(a+1)/n-1}\,ds\,dt =\frac{n}{n-1-a}\|h\|_{L^1}.
\]
The equality for nonnegative $h$ identifies the Hardy-operator norm,
and the asserted constants follow from Theorem~\ref{thm:main} and
Proposition~\ref{prop:sharp}.

(ii) Put $\eta=(a+1)/n-1<0$ and $\lambda=1/p-(a+1)/n>0$. For
$h\geq0$, the Hardy--Littlewood inequality applied to
$s^\eta\chi_{(t,\infty)}(s)$ gives
\[
\HH_{a,b}h(t) \leq t^{(1-\delta)/n}h^{**}(t)+\HH_{a,b}(h^*)(t).
\]
We use the weighted Hardy inequality (see
\cite{KufnerMaligrandaPersson})
\[
\Big\|t^\lambda\int_t^\infty\varphi(s)\,ds\Big\|_{L^r(dt/t)}
\leq\frac1\lambda \big\|t^{\lambda+1}\varphi(t)\big\|_{L^r(dt/t)},
\qquad \lambda>0.
\]
Since $\HH_{a,b}h$ is non-increasing, multiplying the preceding
pointwise estimate by $t^{1/q}$ and taking $L^r(dt/t)$ norms yields
\[
\|\HH_{a,b}h\|_{L^{q,r}}^* \leq\|t^{1/p}h^{**}(t)\|_{L^r(dt/t)}
+\lambda^{-1}\|t^{1/p}h^*(t)\|_{L^r(dt/t)} \leq
(p'+\lambda^{-1})\|h\|_{L^{p,r}}.
\]
Since $q>1$, the starred quasi-norm is equivalent to the $L^{q,r}$
norm. The conclusion follows from Theorem~\ref{thm:main}.
\end{proof}

\begin{remark}\label{rem:lorentz-refined}
For $a=b=0$ we recover the classical Lorentz refinements
\cite{Alvino}, including the sharp constant
$1/((n-1)\omega_n^{1/n})$ when $p=1$. Related results with a
Lebesgue domain appear in \cite[Theorem~8.4]{AlvinoEtAl}. More
generally, for $a\leq0$, the unrestricted and radial optimal
constants coincide by Proposition~\ref{prop:sharp}.
\end{remark}

The optimal-target theorem also identifies the smallest r.i.\ range
for a fixed Lebesgue domain.

\begin{proposition}\label{cor:optimal-lebesgue}
Let $n\geq2$, $b\geq0$, $0\leq\delta=b-a\neq1$, $1\leq p<\infty$,
$a<n/p-1$, and $1/q=1/p-(1-\delta)/n$. If $\delta>1$, assume also
that $q>1$. Then the optimal r.i.\ target for $X=L^p$ in
(\ref{eq:CKN}) is $Y_X=L^{q,p}$, up to equivalent norms. Thus the
choice $r=p$ in Corollary~\ref{cor:lorentz-refined} is optimal among
r.i.\ targets.
\end{proposition}

\begin{proof}
Here $\Xbar'=L^{p'}$, with $p'=\infty$ when $p=1$. Since
$\eta+1/p'=(a+1)/n-1/p<0$, the tail $s^\eta\chi_{(1,\infty)}$
belongs to $L^{p'}$. For $\delta>1$ the additional local condition
$s^{(1-\delta)/n}\chi_{(0,1)}\in L^{p'}$ is equivalent to $q>1$.
Thus Theorem~\ref{thm:optimal} applies. The assumptions also give
$1<q<\infty$.

Monotonicity of $g^*$ gives $\int_0^s g^*(t)t^{-b/n}\,dt
\geq\frac{n}{n-b}s^{1-b/n}g^*(s)$. For the reverse estimate we use
the weighted Hardy inequality (see \cite{KufnerMaligrandaPersson}),
valid when $\eta+1/p'<0$; for $p=1$ it follows directly by
integration. Applied to $\varphi(t)=g^*(t)t^{-b/n}$, the two bounds
become
\[
\frac{n}{n-b}\|s^{(1-\delta)/n}g^*(s)\|_{L^{p'}} \leq \rho(g)
\leq\frac1{|\eta+1/p'|} \|s^{(1-\delta)/n}g^*(s)\|_{L^{p'}}.
\]
Since $(1-\delta)/n=1/q'-1/p'$, this is the $L^{q',p'}$ norm, up to
equivalence. Hence $\Ybar_X'=L^{q',p'}$, and Lorentz duality (see
\cite[Chapter~4]{BS}) gives $Y_X=L^{q,p}$.
\end{proof}

\begin{remark}\label{rem:delta-large}
The classical range $\delta=b-a\leq1$ for Lebesgue CKN inequalities
is not a limitation of the reduction principle. Indeed, suppose that
$b\geq0$, $\delta>1$, $1<p<\infty$, and $1<q<p$, with
$1/q=1/p-(1-\delta)/n$. If $a\geq n/p-1$,
Theorems~\ref{thm:main} and~\ref{thm:optimal} exclude every r.i.\
target for $X=L^p$. If $a<n/p-1$, the optimal
target is $L^{q,p}$ by Proposition~\ref{cor:optimal-lebesgue}; $L^q$
is not admissible because $L^{q,p}\not\subset L^q$. In this range
the boundedness into $L^{q,p}$ is elementary: the weighted Hardy
inequality (the case $\delta=1$) gives $|x|^{-(a+1)}f\in L^p$ and,
since $|x|^{-(\delta-1)}\in L^{n/(\delta-1),\infty}$, O'Neil's
H\"older inequality \cite{ONeil} yields $|x|^{-b}f\in L^{q,p}$. The
content of the reduction principle here is the optimality of this
target among all r.i.\ spaces and the non-existence statement for
$a\geq n/p-1$.

For example, when $n=3$, $p=2$, $a=-2$ and $b=0$, we have $q=6/5$.
Smooth truncations of $f(x)=|x|^{-5/2}(\log|x|)^{-\alpha}$ for
$|x|>e$, with $1/2<\alpha\leq5/6$, have uniformly bounded weighted
gradient norms but unbounded $L^{6/5}$ norms. Weighted Riesz
potentials also yield Lorentz estimates with $q<p$ in certain ranges
\cite{NursultanovTikhonov,LiuYuZhou}; the example above lies outside
their stated power-weight range, since $a=-2<-n/p'=-3/2$.
\end{remark}

\subsection{Critical embeddings and optimality}

At $a=n/p-1$ the strong Lebesgue inequality fails, but a Lorentz
endpoint remains valid with optimal r.i.\ domain and range.

\begin{proposition}\label{cor:critical}
Let $n\geq2$, $\delta\geq0$, $1\leq p<\infty$, and set $a=n/p-1$,
$b=a+\delta$, with $0\leq b<n$. Define $1/q=b/n$, taking $q=\infty$
if $b=0$; thus $1<q\leq\infty$. With the conventions
$L^{\infty,\infty}=L^\infty$ and $q'=1$ if $q=\infty$,
\begin{equation*}
\big\||x|^{-b}f\big\|_{L^{q,\infty}} \leq\frac{q'e^{\Lambda(n,a,b)}}
{n\omega_n^{(1-\delta)/n}} \big\||x|^{-a}\nabla f\big\|_{L^{p,1}},
\qquad f\in C_c^\infty(\R^n).
\end{equation*}
For the norms (\ref{eq:weak-norm}) and (\ref{eq:Lp1}), the optimal
radial constant is $q'/(n\omega_n^{(1-\delta)/n})$; it is also the
optimal general constant when $a\leq0$ or $\delta\geq1$. Among r.i.\
spaces, $L^{p,1}$ is the largest admissible domain for
$L^{q,\infty}$, and $L^{q,\infty}$ is the smallest admissible target
for $L^{p,1}$. For $b=0$ this includes the sharp estimate
\begin{equation}\label{eq:critical-infinity}
\|f\|_{L^\infty} \leq\frac1{n\omega_n^{(1-\delta)/n}}
\big\||x|^\delta\nabla f\big\|_{L^{n/(1-\delta),1}}, \qquad
0\leq\delta<1.
\end{equation}
\end{proposition}

\begin{proof}
The parameter restriction $0\leq b<n$ ensures $1<q\leq\infty$; when
$\delta>1$, it also implies $p>1$ and $1<\delta<1+n/p'$. The Hardy
operator has the form $\HH_{a,b}h(t)=t^{-1/q}\int_t^\infty
h(s)s^{1/p-1}\,ds$. For $h\geq0$ put $I(h)=\int_0^\infty
h(s)s^{1/p-1}\,ds$. Monotonicity gives
\[
\|\HH_{a,b}h\|_{L^{q,\infty}}^*=I(h), \qquad
\|\HH_{a,b}h\|_{L^{q,\infty}}=q'I(h).
\]
For $q<\infty$, the second identity follows from the bound
$\HH_{a,b}h(t)\leq I(h)t^{-1/q}$ and
$\lim_{t\downarrow0}t^{1/q}(\HH_{a,b}h)^{**}(t)=q'I(h)$, by
dominated convergence; for $q=\infty$ it is immediate.
Hardy--Littlewood now yields $I(h)\leq\int_0^\infty
h^*(s)s^{1/p-1}\,ds =\|h\|_{L^{p,1}}$, with equality for decreasing
$h$. Thus $\|\HH_{a,b}\|_{L^{p,1}\to L^{q,\infty}}=q'$. The
inequality and the constants follow from Theorem~\ref{thm:main} and
Proposition~\ref{prop:sharp}. When $b=0$, we have $a=-\delta$,
giving (\ref{eq:critical-infinity}).

If $\HH_{a,b}:\Zbar\to L^{q,\infty}$ is bounded, applying the first
norm identity to $h^*$ shows that
$\|h\|_{L^{p,1}}\preceq\|h\|_{\Zbar}$, hence $Z\hookrightarrow
L^{p,1}$. For the target, suppose that $\HH_{a,b}:L^{p,1}\to\Zbar$
is bounded. Taking $h_R=R^{-1/p}\chi_{(R,2R)}$ gives
$\|h_R\|_{L^{p,1}}=p$ and, for $0<t<R$,
\[
\HH_{a,b}h_R(t) =p(2^{1/p}-1)t^{-1/q}.
\]
The lattice and Fatou properties imply $t^{-1/q}\in\Zbar$ (the
constant function $1$ when $q=\infty$). Since every $k\in
L^{q,\infty}$ satisfies $k^*(t)\leq\|k\|_{L^{q,\infty}}^*t^{-1/q}$,
it follows that $L^{q,\infty}\hookrightarrow Z$.
Theorem~\ref{thm:main} transfers both optimality assertions to CKN.
\end{proof}

\begin{remark}\label{rem:critical}
For $b>0$, $bq=n$ and the measure $|x|^{-bq}dx=|x|^{-n}dx$ has
infinite mass near the origin. Thus, for smooth $f$ with
$f(0)\neq0$, the weighted weak-$L^q$ norm is infinite, whereas the
multiplier norm $\||x|^{-b}f\|_{L^{q,\infty}}$ is finite.

Moreover, at $a=n/p-1$ and $p>1$ we have $\eta=-1/p'$, so
$s^\eta\chi_{(1,\infty)}\notin L^{p'}$.
Theorem~\ref{thm:optimal}(i) therefore rules out every r.i.\ target
for $X=L^p$: the critical domain $L^{p,1}$ cannot be replaced by
$L^p$. The general problem of determining optimal domains for
arbitrary prescribed r.i.\ targets belongs to the weighted
Hardy-operator theory developed in \cite{Mihula}.

The case $p=n$, $a=0$, $b=1$ deserves an explicit comment. The
preceding paragraph shows that, for $n\geq2$, the weak-type Hardy
inequality
\begin{equation}\label{eq:weak-hardy-n}
\Big\|\frac{f}{|x|}\Big\|_{L^{n,\infty}(\R^n)}\leq C\|\nabla
f\|_{L^n(\R^n)}, \qquad f\in C_c^\infty(\R^n),
\end{equation}
cannot hold. This can also be seen directly: for $R>1$ let
$f_R(x)=\log R$ if $|x|\leq1$, $f_R(x)=\log(R/|x|)$ if $1<|x|<R$,
and $f_R(x)=0$ otherwise. Then $\|\nabla f_R\|_{L^n}^n=n\omega_n\log
R$, whereas $|x|^{-1}f_R\geq(\log R)|x|^{-1}\chi_{\{|x|<1\}}$ and
$\||x|^{-1}\chi_{\{|x|<1\}}\|_{L^{n,\infty}}^*=\omega_n^{1/n}$.
Hence the quotient of the two sides of (\ref{eq:weak-hardy-n}) is at
least a constant multiple of $(\log R)^{1-1/n}$, and a standard
regularization of $f_R$ gives the same conclusion within
$C_c^\infty(\R^n)$. Thus the weak-type endpoint requires the smaller
domain $L^{n,1}$, as established in Proposition~\ref{cor:critical}.
\end{remark}

\subsection{Weighted Hardy inequalities}

For $\delta=1$ the operator becomes
$\HH_{b-1,b}h(t)=t^{-b/n}\int_t^\infty h(s)s^{b/n-1}\,ds$. An
application to arbitrary r.i.\ spaces is conveniently expressed in
terms of dilations. Write $E_s h(t)=h(t/s)$,
$h_X(s)=\|E_s\|_{\Xbar\to\Xbar}$, and let
$\underline\alpha_X=\lim_{s\to0^+}\log h_X(s)/\log s$ be the lower
Boyd index of $X$ (see \cite[Chapter~3, Section~5]{BS}).

\begin{proposition}\label{cor:hardy-ri}
Let $n\geq2$, $b\geq0$, and let $X$ be an r.i.\ space with
$\underline\alpha_X>b/n$. Then
\begin{equation}\label{eq:hardy-ri}
\big\||x|^{-b}f\big\|_X \leq\frac1n\left(\int_1^\infty
h_X(1/u)u^{b/n-1}\,du\right) \big\||x|^{1-b}\nabla f\big\|_X, \qquad
f\in C_c^\infty(\R^n),
\end{equation}
and the integral is finite. For $X=L^{p,r}$, with $(p,r)=(1,1)$ or
$1<p<\infty$, $1\leq r\leq\infty$, and $b<n/p$, equipped with an
r.i.\ norm for which $\|E_s h\|=s^{1/p}\|h\|$, the optimal constant
is $p/(n-bp)$.
\end{proposition}

\begin{proof}
The substitution $s=tu$ gives $\HH_{b-1,b}h=\int_1^\infty
E_{1/u}h\,u^{b/n-1}\,du$. Minkowski's integral inequality in $\Xbar$
therefore yields
\[
\|\HH_{b-1,b}h\|_{\Xbar} \leq\left(\int_1^\infty
h_X(1/u)u^{b/n-1}\,du\right) \|h\|_{\Xbar}.
\]
The definition of $\underline\alpha_X>b/n$ makes the integral
finite. Since $\delta=1$, the equality of the radial and general
optimal constants in Proposition~\ref{prop:sharp} gives
(\ref{eq:hardy-ri}).

For the Lorentz spaces in the statement, $h_X(s)=s^{1/p}$ and the
integral equals $np/(n-bp)$. To see that the operator norm attains
this upper bound asymptotically, use
$h_\varepsilon(t)=t^{-1/p-\varepsilon}\chi_{(1,\infty)}(t)$. These
functions belong to $L^{p,r}$; for $t\geq1$,
$\HH_{b-1,b}h_\varepsilon(t)
=(1/p+\varepsilon-b/n)^{-1}h_\varepsilon(t)$. Letting
$\varepsilon\downarrow0$ proves sharpness.
\end{proof}

In particular, $X=L^p$ gives the classical sharp weighted Hardy
constant $p/(n-bp)$, including $b=1$ when $p<n$.

\begin{remark}\label{rem:hardy-direct}
Inequality (\ref{eq:hardy-ri}) also admits a direct proof, which
does not use rearrangements of $f$ and is valid for every real $b$
for which the integral is finite. Indeed, for $f\in
C_c^\infty(\R^n)$ and $x\neq0$, integration along the ray through
$x$ gives $|f(x)|\leq\int_1^\infty|\nabla f(\lambda
x)|\,|x|\,d\lambda$, that is,
\[
|x|^{-b}|f(x)|\leq\int_1^\infty\lambda^{b-1}G(\lambda x)\,d\lambda,
\qquad G(y)=|y|^{1-b}|\nabla f(y)|.
\]
Since $\bigl(G(\lambda\,\cdot)\bigr)^*(t)=G^*(\lambda^nt)$,
Minkowski's integral inequality and the substitution $u=\lambda^n$
yield
\[
\big\||x|^{-b}f\big\|_X\leq\int_1^\infty\lambda^{b-1}
h_X(\lambda^{-n})\,d\lambda\,\|G\|_X
=\frac1n\left(\int_1^\infty h_X(1/u)u^{b/n-1}\,du\right)\|G\|_X .
\]
What Theorem~\ref{thm:main} and Proposition~\ref{prop:sharp} add to
this elementary argument is the converse: for $b\geq0$, the weighted
Hardy inequality holds in $X$ if and only if $\HH_{b-1,b}$ is
bounded on $\Xbar$, with optimal constant
$\|\HH_{b-1,b}\|_{\Xbar\to\Xbar}/n$.
\end{remark}

\section*{Statements and Declarations}

\noindent\textbf{Funding.} The first author is partially supported
by Grants PID2024-160507NB-I00 and PID2024-155917NB-I00, funded by
MCIN/AEI/10.13039/501100011033.

The second author is partially supported by Grant
PID2024-160507NB-I00, funded by MCIN/AEI/10.13039/501100011033.

\medskip
\noindent\textbf{Competing interests.} The author declares no
competing financial or non-financial interests.

\medskip
\noindent\textbf{Data availability.} No datasets were generated or
analysed during the current study.

\medskip
\noindent\textbf{Use of generative AI.}

The main reduction principle, based on relative rearrangements and
an integrating factor, together with its original proof, was
developed by the authors without the use of generative AI.

During the revision process, ChatGPT (OpenAI) and Claude (Anthropic)
were used to assist with bibliographic searches, improve the
exposition, and discuss and critically examine mathematical
arguments. AI-assisted suggestions contributed to specific
refinements, including the sharper constant in
Lemma~\ref{lem:pointwise} and the perturbative argument in
Proposition~\ref{prop:perturbative}.

All AI-assisted mathematical contributions were independently
verified by the authors, who take full responsibility for the entire
content of the paper.

\end{document}